\documentclass[12pt]{article}
\usepackage{amsfonts}
\usepackage{full page,
amssymb, amscd, amsmath, graphicx, 
}
\allowdisplaybreaks

 \title{{\bf Associative algebras and the representation theory of 
grading-restricted vertex algebras}}
 \author{Yi-Zhi Huang}
    \date{}
    \begin{document}
    \bibliographystyle{alpha}
    \maketitle
\newtheorem{thm}{Theorem}[section]
\newtheorem{defn}[thm]{Definition}
\newtheorem{prop}[thm]{Proposition}
\newtheorem{cor}[thm]{Corollary}
\newtheorem{lemma}[thm]{Lemma}
\newtheorem{rema}[thm]{Remark}
\newtheorem{app}[thm]{Application}
\newtheorem{prob}[thm]{Problem}
\newtheorem{conv}[thm]{Convention}
\newtheorem{conj}[thm]{Conjecture}
\newtheorem{cond}[thm]{Condition}
    \newtheorem{exam}[thm]{Example}
\newtheorem{assum}[thm]{Assumption}
     \newtheorem{nota}[thm]{Notation}
\newcommand{\halmos}{\rule{1ex}{1.4ex}}
\newcommand{\pfbox}{\hspace*{\fill}\mbox{$\halmos$}}
\newcommand{\overarc}[1]{\,\arc{\!{#1}}}
\newcommand{\nn}{\nonumber \\}
\newcommand{\overbar}[1]{\,\overline{\!{#1}}}
 \newcommand{\res}{\mbox{\rm Res}}
 \newcommand{\ord}{\mbox{\rm ord}}
\renewcommand{\hom}{\mbox{\rm Hom}}
\newcommand{\edo}{\mbox{\rm End}\ }
 \newcommand{\pf}{{\it Proof.}\hspace{2ex}}
 \newcommand{\epf}{\hspace*{\fill}\mbox{$\halmos$}}
 \newcommand{\epfv}{\hspace*{\fill}\mbox{$\halmos$}\vspace{1em}}
 \newcommand{\epfe}{\hspace{2em}\halmos}
\newcommand{\nord}{\mbox{\scriptsize ${\circ\atop\circ}$}}
\newcommand{\wt}{\mbox{\rm wt}\ }
\newcommand{\swt}{\mbox{\rm {\scriptsize wt}}\ }
\newcommand{\lwt}{\mbox{\rm wt}^{L}\;}
\newcommand{\rwt}{\mbox{\rm wt}^{R}\;}
\newcommand{\slwt}{\mbox{\rm {\scriptsize wt}}^{L}\,}
\newcommand{\srwt}{\mbox{\rm {\scriptsize wt}}^{R}\,}
\newcommand{\clr}{\mbox{\rm clr}\ }
\newcommand{\tr}{\mbox{\rm Tr}}
\renewcommand{\H}{\mathbb{H}}
\newcommand{\C}{\mathbb{C}}
\newcommand{\Z}{\mathbb{Z}}
\newcommand{\R}{\mathbb{R}}
\newcommand{\Q}{\mathbb{Q}}
\newcommand{\N}{\mathbb{N}}
\newcommand{\CN}{\mathcal{N}}
\newcommand{\F}{\mathcal{F}}
\newcommand{\I}{\mathcal{I}}
\newcommand{\V}{\mathcal{V}}
\newcommand{\one}{\mathbf{1}}
\newcommand{\BY}{\mathbb{Y}}
\newcommand{\ds}{\displaystyle}
\newcommand{\g}{\textsl{g}}

        \newcommand{\ba}{\begin{array}}
        \newcommand{\ea}{\end{array}}
        \newcommand{\be}{\begin{equation}}
        \newcommand{\ee}{\end{equation}}
        \newcommand{\bea}{\begin{eqnarray}}
        \newcommand{\eea}{\end{eqnarray}}
         \newcommand{\lbar}{\bigg\vert}
      \newcommand{\mbar}{\big\vert}
        \newcommand{\p}{\partial}
        \newcommand{\dps}{\displaystyle}
        \newcommand{\bra}{\langle}
        \newcommand{\ket}{\rangle}

        \newcommand{\ob}{{\rm ob}\,}
        \renewcommand{\hom}{{\rm Hom}}

\newcommand{\A}{\mathcal{A}}
\newcommand{\Y}{\mathcal{Y}}
\newcommand{\End}{\mathrm{End}}
 \newcommand{\rad}{\mbox{\rm rad}}

\newcommand{\dlt}[3]{#1 ^{-1}\delta \bigg( \frac{#2 #3 }{#1 }\bigg) }

\newcommand{\dlti}[3]{#1 \delta \bigg( \frac{#2 #3 }{#1 ^{-1}}\bigg) }

 \makeatletter
\newlength{\@pxlwd} \newlength{\@rulewd} \newlength{\@pxlht}
\catcode`.=\active \catcode`B=\active \catcode`:=\active
\catcode`|=\active
\def\sprite#1(#2,#3)[#4,#5]{
   \edef\@sprbox{\expandafter\@cdr\string#1\@nil @box}
   \expandafter\newsavebox\csname\@sprbox\endcsname
   \edef#1{\expandafter\usebox\csname\@sprbox\endcsname}
   \expandafter\setbox\csname\@sprbox\endcsname =\hbox\bgroup
   \vbox\bgroup
  \catcode`.=\active\catcode`B=\active\catcode`:=\active\catcode`|=\active
      \@pxlwd=#4 \divide\@pxlwd by #3 \@rulewd=\@pxlwd
      \@pxlht=#5 \divide\@pxlht by #2
      \def .{\hskip \@pxlwd \ignorespaces}
      \def B{\@ifnextchar B{\advance\@rulewd by \@pxlwd}{\vrule
         height \@pxlht width \@rulewd depth 0 pt \@rulewd=\@pxlwd}}
      \def :{\hbox\bgroup\vrule height \@pxlht width 0pt depth
0pt\ignorespaces}
      \def |{\vrule height \@pxlht width 0pt depth 0pt\egroup
         \prevdepth= -1000 pt}
   }
\def\endsprite{\egroup\egroup}
\catcode`.=12 \catcode`B=11 \catcode`:=12 \catcode`|=12\relax
\makeatother

\def\hboxtr{\FormOfHboxtr} 
\sprite{\FormOfHboxtr}(25,25)[0.5 em, 1.2 ex] 

:BBBBBBBBBBBBBBBBBBBBBBBBB | :BB......................B |
:B.B.....................B | :B..B....................B |
:B...B...................B | :B....B..................B |
:B.....B.................B | :B......B................B |
:B.......B...............B | :B........B..............B |
:B.........B.............B | :B..........B............B |
:B...........B...........B | :B............B..........B |
:B.............B.........B | :B..............B........B |
:B...............B.......B | :B................B......B |
:B.................B.....B | :B..................B....B |
:B...................B...B | :B....................B..B |
:B.....................B.B | :B......................BB |
:BBBBBBBBBBBBBBBBBBBBBBBBB |

\endsprite
\def\shboxtr{\FormOfShboxtr} 
\sprite{\FormOfShboxtr}(25,25)[0.3 em, 0.72 ex] 

:BBBBBBBBBBBBBBBBBBBBBBBBB | :BB......................B |
:B.B.....................B | :B..B....................B |
:B...B...................B | :B....B..................B |
:B.....B.................B | :B......B................B |
:B.......B...............B | :B........B..............B |
:B.........B.............B | :B..........B............B |
:B...........B...........B | :B............B..........B |
:B.............B.........B | :B..............B........B |
:B...............B.......B | :B................B......B |
:B.................B.....B | :B..................B....B |
:B...................B...B | :B....................B..B |
:B.....................B.B | :B......................BB |
:BBBBBBBBBBBBBBBBBBBBBBBBB |

\endsprite

\makeatletter
\DeclareFontFamily{U}{tipa}{}
\DeclareFontShape{U}{tipa}{m}{n}{<->tipa10}{}
\newcommand{\arc@char}{{\usefont{U}{tipa}{m}{n}\symbol{62}}}%

\newcommand{\arc}[1]{\mathpalette\arc@arc{#1}}

\newcommand{\arc@arc}[2]{%
  \sbox0{$\m@th#1#2$}%
  \vbox{
    \hbox{\resizebox{\wd0}{\height}{\arc@char}}
    \nointerlineskip
    \box0
  }%
}
\makeatother

\date{}
\maketitle

\begin{abstract}
We introduce an associative algebra $A^{\infty}(V)$ using infinite matrices with entries in 
a grading-restricted  vertex algebra $V$ such that the associated graded 
space $Gr(W)=\coprod_{n\in \N}Gr_{n}(W)$ of a filtration of a lower-bounded 
generalized $V$-module $W$ is an $A^{\infty}(V)$-module satisfying additional 
properties (called a nondegenerate 
graded $A^{\infty}(V)$-module). We prove that a lower-bounded generalized $V$-module $W$
is irreducible or completely reducible if and only if the nondegenerate graded $A^{\infty}(V)$-module
$Gr(W)$ is irreducible or completely reducible, respectively. We also prove that 
the set of equivalence classes of the lower-bounded generalized $V$-modules are in bijection with the set of the
equivalence classes of nondegenerate graded $A^{\infty}(V)$-modules. 
For $N\in \N$, there is a subalgebra $A^{N}(V)$ of $A^{\infty}(V)$ 
such that the subspace $Gr^{N}(W)=\coprod_{n=0}^{N}Gr_{n}(W)$
of $Gr(W)$ is an $A^{N}(V)$-module satisfying additional properties (called a 
nondegenerate graded $A^{N}(V)$-module).
We prove that $A^{N}(V)$  are finite dimensional 
when $V$ is of positive energy  (CFT type) and $C_{2}$-cofinite. We prove that the set of 
the equivalence classes of
lower-bounded generalized $V$-modules is in bijection with the set of the equivalence 
classes of nondegenerate graded $A^{N}(V)$-modules. 
In the case that $V$ is a M\"{o}bius vertex algebra and 
the differences between the real parts of the lowest weights of the 
irreducible lower-bounded generalized $V$-modules are less than or equal to $N\in \N$,  we prove that
 a lower-bounded generalized $V$-module $W$ of finite length
is irreducible or completely reducible if and only if the nondegenerate graded $A^{N}(V)$-module
$Gr^{N}(W)$ is irreducible or completely reducible, respectively. 
\end{abstract}



\renewcommand{\theequation}{\thesection.\arabic{equation}}
\renewcommand{\thethm}{\thesection.\arabic{thm}}
\setcounter{equation}{0} \setcounter{thm}{0} 

\section{Introduction}

In the representation theory of Lie algebras, the universal enveloping 
algebra of a Lie algebra plays a crucial role because the module categories 
for a Lie algebra and for its universal enveloping algebra are isomorphic. 
For a vertex operator algebra, there is also a universal enveloping algebra 
introduced by Frenkel and Zhu \cite{FZ} such that the module categories for 
these algebras are isomorphic. Unfortunately, the universal enveloping 
algebra of a vertex operator algebra  is not very useful since it involves certain infinite sums in a 
suitable topological completion of the tensor algebra of the algebra.
On the other hand,  the classes of modules that we are interested in 
the representation theory of vertex operator algebras and two-dimensional 
conformal field theory do not involve such infinite sums since
the vertex operators on these modules are lower truncated when acting 
on elements of these modules.

Instead, in the representation theory of vertex operator algebras, we 
have the Zhu algebra $A(V)$ introduced by Zhu in \cite{Z} and it's generalizations $A_{n}(V)$ for 
$n\in N$ by Dong, Li and Mason in \cite{DLM} for a vertex operator algebra $V$. These 
algebras can be used to classify irreducible modules for the vertex operator 
algebra and to study problems related to different types of modules. But
compared with the universal enveloping algebra of a Lie algebra, the role 
of these associative algebras played in the representation theory of vertex operator algebras 
is quite limited. For example, the module for one of these associative algebras obtained from a 
suitable $V$-module in general do not tell us whether the original $V$-module is irreducible 
or completely reducible. 

In the present paper, we introduce an associative algebra $A^{\infty}(V)$ 
using infinite matrices with entries in a grading-restricted vertex algebra $V$. 
The associated graded space $Gr(W)=\coprod_{n\in \N}Gr_{n}(W)$ of a filtration of 
a lower-bounded generalized $V$-module $W$ is an $A^{\infty}(V)$-module with an $\N$-grading and
some operators having suitable properties (called a nondegenerate graded $A^{\infty}(V)$-module). 
In fact, the algebra $A^{\infty}(V)$ is defined using the associated graded spaces of
all  lower-bounded generalized $V$-modules. 
We prove that a lower-bounded generalized $V$-module $W$
is irreducible or completely reducible if and only if the nondegenerate graded $A^{\infty}(V)$-module
$Gr(W)$ is irreducible or completely reducible, respectively. We also prove that 
the set of the equivalence classes of irreducible lower-bounded generalized $V$-modules
is in bijection with the set of the equivalence classes of irreducible nondegenerate graded $A^{\infty}(V)$-modules. 

We show that $A(V)$ in \cite{Z} and $A_{n}(V)$ \cite{DLM} mentioned above are isomorphic to
very special subalgebras of 
$A^{\infty}(V)$. This fact gives a conceptual explanation of the role that these associative
algebras played in the representation theory of vertex operator algebras.

We then introduce new subalgebras $A^{N}(V)$ of 
$A^{\infty}(V)$ for $N\in \N$. These subalgebras can also be obtained using 
finite matrices with entries in $V$. In the case that $V$ is of positive energy (or CFT type) and 
$C_{2}$-cofinite, we prove that $A^{N}(V)$  are finite dimensional. 
The subspace $Gr^{N}(W)=\coprod_{n=0}^{N}Gr_{n}(W)$ of 
$Gr(W)$ of a lower-bounded generalized $V$-module $W$ is an $A^{N}(V)$-module with
some operators having suitable properties (called a nondegenerate graded $A^{N}(V)$-module). 
Wer prove that if a lower-bounded generalized $V$-module $W$ is irreducible or completely reducible, 
then the nondegenerate graded $A^{N}(V)$-module $Gr^{N}(W)$ is irreducible or completely reducible, respectively.  
We also prove that the set of the
equivalence classes of irreducible lower-bounded generalized $V$-modules
is in bijection with the set of the equivalence classes of 
irreducible nondegenerate graded $A^{N}(V)$-modules.  

In the case that $V$ is a M\"{o}bius vertex 
algebra so that a lowest weight of a lower-bounded generalized $V$-module is defined, 
under the assumption that the differences between the 
real parts of the lowest weights of the irreducible lower-bounded generalized 
$V$-modules are less than or equal to $N\in \N$, we prove that 
 a lower-bounded generalized $V$-module $W$ of finite length
is irreducible or completely reducible if and only if the nondegenerate graded $A^{N}(V)$-module
$Gr^{N}(W)$ is irreducible or completely reducible, respectively. 
When $A^{N}(V)$ for all $N\in \N$ are finite dimensional (for example, when 
$V$ is of positive energy (or CFT type) and 
$C_{2}$-cofinite), we prove that 
an irreducible lower-bounded generalized $V$-module is an ordinary $V$-module and 
thus every lower-bounded generalized $V$-module of finite length is grading-restricted. 
In this case, under the assumptions
above on $V$, lowest weights and $N$, a lower-bounded generalized $V$-module 
$W$ of finite length or a grading-restricted 
generalized $V$-module $W$ is a direct sum of irreducible ordinary $V$-modules
 if and only if  the nondegenerate graded $A^{N}(V)$-module
$Gr^{N}(W)$ is completely reducible. 

Many of the main results mentioned above need the construction of 
universal lower-bounded generalized $V$-modules in \cite{H-const-twisted-mod} and 
some results from \cite{H-exist-twisted-mod}. 

The category of lower-bounded generalized $V$-modules and the category of nondegenerate
graded $A^{\infty}(V)$-modules are not equivalent because of morphisms, but they are ``almost'' equivalent. 
We shall study the relations between these categories, the category of lower-bounded 
generalized $V$-modules of finite lengths  and the categories of graded 
$A^{N}(V)$-modules for $N\in \N$
in another paper.

This paper is organized as follows: In the next section, we introduce the associative algebra 
$A^{\infty}(V)$ associated to a grading-restricted vertex algebra $V$
 and prove that the associated graded space $Gr(W)$ of a filtration of 
a lower-bounded generalized $V$-module $W$ is an $A^{\infty}(V)$-module. 
In Section 3, we introduce graded $A^{\infty}(V)$-modules and 
prove the results mentioned above on the relations between 
 lower-bounded generalized $V$-modules and nondegenerate graded $A^{\infty}(V)$-modules.
We show that the Zhu 
algebra and their generalizations by Dong, Li and Mason are isomorphic to 
subalgebras of $A^{\infty}(V)$ in
Subsection 4.1 and introduce the new subalgebras $A^{N}(V)$ of $A^{\infty}(V)$ for $N\in \N$ in Subsection 4.2. 
We also prove in Subsection 4.2 that when $V$ is of positive energy and 
$C_{2}$-cofinite, $A^{N}(V)$ for $N\in \N$ are finite dimensional. 
In Section 5, we introduce graded $A^{N}(V)$-modules and 
prove the results mentioned above on the relations between lower-bounded generalized 
$V$-modules, lower-bounded generalized 
$V$-modules of finite lengths and nondegenerate graded $A^{N}(V)$-modules. 

\paragraph{Acknowledgment} I am very grateful to Darlayne Addabbo
for noticing that an early version of Theorem \ref{DLM-isom} (which has been 
corrected now) does not hold
for the Heisenberg vertex operator algebra.

\setcounter{equation}{0} \setcounter{thm}{0} 
\section{Associative algebra $A^{\infty}(V)$ and modules}

In this paper, we fix a grading-restricted vertex algebra $V$. Most  of the constructions and 
results work and hold for more general vertex algebras, for example, lower-bounded vertex algebras or 
superalgebras. The constructions and results certainly work and hold for a M\"{o}bius vertex algebra or a
vertex operator algebra. For some results in Section 5, 
we shall assume that $V$ is a M\"{o}bius vertex algebra.

Let $U^{\infty}(\C)$ be the space of all  column-finite infinite matrices 
with entries in $\C$, but doubly indexed by $\N$ instead of $\Z_{+}$. 
In other words, $U^{\infty}(\C)$ is the space of all infinite matrices
of the form $[a_{kl}]$ for $a_{kl}\in \C$, $k, l\in \N$ such that for each fixed $l\in \N$, 
there are only finitely many nonzero $a_{kl}$. Let $I^{\infty}=[\delta_{kl}]$ be the infinite identity matrix.
Then $U^{\infty}(\C)$ is in fact an 
associative algebra with the identity $I^{\infty}$. 
The space $U^{\infty}(\C)$ has a set of linearly independent elements of the form
$E_{kl}$ for $k, l\in \N$ with the entry in the $k$-th row and $l$-th column equal to $1$ and 
all the other entries equal to $0$. 
These infinite matrices do not form a basis of $U^{\infty}(\C)$ but  form a basis of the subspace
$U^{\infty}_{0}(\C)$ of 
$U^{\infty}(\C)$ consisting of finitary matrices (matrices with only finitely many nonzero entries). In particular, 
$$U^{\infty}_{0}(\C)=\coprod_{k, l\in \N}\C E_{kl}.$$
Moreover, 
$$U^{\infty}(\C)\subset \prod_{k, l\in \N}\C E_{kl}, $$
where $\prod_{k, l\in \N}\C E_{kl}$ is the algebraic completion of $U^{\infty}_{0}(\C)$ viewed as a graded space. 
Though elements of $U^{\infty}(\C)$ are infinite linear combinations of 
$E_{kl}$ for $k, l\in \N$, any  binary product on $U^{\infty}(\C)$  satisfying the distribution axioms 
is still determined completely by the product of 
$E_{kl}$ for $k, l\in \N$. For example, for the usual  matrix product, 
we know that $E_{kl}E_{mn}=\delta_{lm}E_{kn}$ for $k, l, m, n\in \N$. 
Let $A=\sum_{k, l\in \N}a_{kl}E_{kl}$ and $B=\sum_{k, l\in \N}b_{kl}E_{kl}$,
where $a_{kl}, b_{kl}\in \C$ for $k, l\in \N$.
Then 
\begin{align*}
AB=\left(\sum_{k, l\in \N}a_{kl}E_{kl}\right)\left(\sum_{m, n\in \N}b_{mn}E_{mn}\right)
=\sum_{k, n\in \N}\left(\sum_{m\in \N}a_{km}b_{mn}\right)E_{kn}.
\end{align*}
So even though $E_{kl}$ for $k, l\in \N$ do not form a basis of $U^{\infty}(\C)$,
all the properties of the associative algebra structure on $U^{\infty}(\C)$ can still be derived from 
the properties these matrices. Thus we can study $U^{\infty}(\C)$ using $E_{kl}$ for $k, l\in \N$, $k\le l$.
Also what we are mainly interested is the subalgebra $\C I^{\infty}\oplus U^{\infty}_{0}(\C)$ of 
$U^{\infty}(\C)$. This subalgebra has a basis $\{I^{\infty}\}\cup \{E_{kl}\}_{k, l\in \N}$. 

Let $U^{\infty}(V)=V\otimes U^{\infty}(\C)$. Then $U^{\infty}(V)$ is  
the space of  column-finite infinite
matrices with entries in $V$, but doubly indexed by $\N$ instead of $\Z_{+}$. 
Elements of  $U^{\infty}(V)$ are of the form 
$\mathfrak{v}=[v_{kl}]$
for $v_{kl}\in V$, $k, l\in \N$  such that for each fixed  $l\in \N$, 
there are only finitely many nonzero $v_{kl}$. 
Let $U^{\infty}_{0}(V)$ be the subspace of $U^{\infty}(V)$ spanned by elements of the form 
$v\otimes E_{kl}$ for $v\in V$ and $k, l\in \N$.   Then 
$$U^{\infty}_{0}(V)=\coprod_{k, l\in \N}V\otimes \C E_{kl}$$
and 
$$U^{\infty}(V)\subset \prod_{k, l\in \N}V\otimes \C E_{kl}.$$
We shall denote 
$v\otimes E_{kl}$ simply by $[v]_{kl}$. Then elements of $U^{\infty}(\C)$ can all be written as 
$$\sum_{k, l\in \N}[v_{kl}]_{kl}$$
for $v_{kl}\in V$, $k, l\in \N$.
As in the case of $U^{\infty}(\C)$,
we can study any binary product on $U^{\infty}(V)$  satisfying the distribution axioms 
using $[v]_{kl}$ for $v\in V$ and $k, l\in \N$. 
We are also mainly interested in the subspace $V\otimes I^{\infty}\oplus U^{\infty}_{0}(V)$ 
of $U^{\infty}(V)$. This subspace is spanned by elements of the form $v\otimes I^{\infty}$ and 
$[v]_{kl}$ for $v\in V$ and $k, l\in \N$. Because of this reason, though we might give definitions
of products and related notions using general elements of $U^{\infty}(V)$, we shall study them 
using only $[v]_{kl}$ for $v\in V$ and $k, l\in \N$. 

We also need some particular formal series and polynomials. 
In this paper, we shall use the convention that 
a complex power or the integral power of the  logarithm of an ordered linear combination of
formal variables and a complex number, 
always means its expansion in nonnegative powers of the formal variables 
or the complex number that are not the first one in the ordered linear combination.
For example, $(x+1)^{-k-1}$ for $k\in \N$ and $(1+x)^{n}$ for $n\in \N$ 
mean the expansions in nonnegative powers of $1$ and in nonnegative powers of $x$,
respectively. For $k, n,  l\in \N$, we have 
\begin{align}\label{exp-x+1-k-1}
(x+1)^{-k+n-l-1}&=\sum_{m\in \N}\binom{-k+n-l-1}{m}x^{-k+n-l-m-1}\nn
&=T_{k+l+1}((x+1)^{-k+n-l-1})+R_{k+l+1}((x+1)^{-k+n-l-1}),
\end{align}
where 
$$T_{k+l+1}((x+1)^{-k+n-l-1})=\sum_{m=0}^{n}\binom{-k+n-l-1}{m}x^{-k+n-l-m-1}$$
is the Taylor polynomial in $x^{-1}$ of order $k+l+1$ of $(x+1)^{-k+n-l-1}$ and 
$$R_{k+l+1}((x+1)^{-k+n-l-1})=\sum_{m\in n+1+\N}\binom{-k+n-l-1}{m}x^{-k+n-l-m-1}$$
is the remainder of order $k+l+1$. 

We define a product $\diamond$ on  $U^{\infty}(V)$ by
$$\mathfrak{u}\diamond \mathfrak{v}=[(\mathfrak{u}\diamond\mathfrak{v})_{kl}]$$
for $\mathfrak{u}=[u_{kl}], \mathfrak{v}=[v_{kl}]\in U^{\infty}(V)$, where 
\begin{align}\label{defn-diamond}
(\mathfrak{u}\diamond\mathfrak{v})_{kl}&=\sum_{n=k}^{l}\res_{x}T_{k+l+1}((x+1)^{-k+n-l-1})
(1+x)^{l}Y_{V}((1+x)^{L_{V}(0)}u_{kn}, x)v_{nl}\nn
&=\sum_{n=k}^{n}\sum_{m=0}^{l}
\binom{-k+n-l-1}{m}\res_{x}
x^{-k+n-l-m-1}(1+x)^{l}Y_{V}((1+x)^{L_{V}(0)}u_{kn}, x)v_{nl}
\end{align}
for $k, l\in \N$. Then $U^{\infty}(V)$ equipped with $\diamond$ is an algebra but in general is not 
even associative. 
Let $O^{\infty}(V)$ be the subspace of $U^{\infty}(V)$ spanned by elements of the form 
$$\left[\sum_{n=k}^{l}\res_{x}
x^{-k-l-p-2}(1+x)^{l}Y_{V}((1+x)^{L_{V}(0)}u_{kn}, x)v_{nl}\right]$$ 
for $\mathfrak{u}=[u_{kl}], \mathfrak{v}=[v_{kl}]\in U^{\infty}(V)$, 
$p\in \N$ and elements of the form 
$$[(L_{V}(-1)+L_{V}(0)+l-k)v_{kl}]$$
 for $\mathfrak{v}=[v_{kl}]\in U^{\infty}(V)$. 

The product $\diamond$ on $U^{\infty}(V)$ looks complicated. But as we mentioned above, 
though $[v]_{kl}$ for $v\in V$ and $k, l\in \N$ does span $U^{\infty}(V)$, 
their infinite linear combinations give all the elements of $U^{\infty}(V)$ and 
$U^{\infty}(V)$ can be studied using these elements. In particular, the product $\diamond$ can be 
studied using these elements. So instead of working with
arbitrary matrices in $U^{\infty}(V)$, we use 
$[v]_{kl}$ for $v\in V$ and $k, l\in \N$ to write down $\diamond$.
For $u, v\in V$ and $k, m, n, l\in \N$, by definition, 
$$[u]_{km}\diamond[v]_{nl}=0$$
when $m\ne n$ and 
\begin{align}\label{defn-diamond-1}
[u]_{kn}\diamond[v]_{nl}&=
\res_{x}T_{k+l+1}((x+1)^{-k+n-l-1})(1+x)^{l}\left[Y_{V}((1+x)^{L_{V}(0)}u, x)v\right]_{kl}\nn
&=\sum_{m=0}^{n}
\binom{-k+n-l-1}{m}\res_{x}
x^{-k+n-l-m-1}(1+x)^{l}\left[Y_{V}((1+x)^{L_{V}(0)}u, x)v\right]_{kl}.
\end{align}
Since $[u]_{km}\diamond[v]_{nl}=0$ when $m\ne n$, we need only consider 
$[u]_{kn}\diamond[v]_{nl}$ for $u, v\in V$ and $k, n, l\in \N$.
By taking $\mathfrak{u}=[u]_{kn}$ and $\mathfrak{v}=[v]_{nl}$,
we see also that the subspace
$O^{\infty}(V)$ is spanned by infinite linear combinations of elements of the form 
$$ \res_{x}
x^{-k-l-p-2}(1+x)^{l}[Y_{V}((1+x)^{L_{V}(0)}u, x)v]_{kl}$$
for $u, v\in V$,  $k, l, p\in \N$ and 
elements of the form 
$$[(L_{V}(-1)+L_{V}(0)+l-k)v]_{kl}$$
for $v\in V$ and $k, l\in \N$, with each pair $(k, l)$ appearing in the 
linear combinations only finitely many times. 

Let $\one^{\infty}$ be the element of $U^{\infty}(V)$ with diagonal entries being $\one \in V$ and
all the other entries being $0$. Then $\one^{\infty}=\one\otimes I^{\infty}$.

We shall take a quotient of $U^{\infty}(V)$ such that the quotient with the product induced 
from $\diamond$ is an associative algebra
and such that  the associated graded space of a filtration of
every lower-bounded generalized $V$-module is a module for this associative algebra. 
To do this, we need to first give an action of the (nonassociative) algebra $U^{\infty}(V)$ with the product $\diamond$
on a lower-bounded generalized $V$-module. 

We briefly recall the notion of lower-bounded generalized $V$-module. We refer the reader to Definition 1.2
in \cite{H-cofiniteness}, where a lower-bounded generalized $V$-module is called a lower-truncated 
generalized $V$-module.  Definition 1.2
in \cite{H-cofiniteness} is for a vertex operator algebra $V$ but the definition applies also 
to a grading-restricted vertex algebra except that we have to 
require the existence of operators $L_{W}(0)$ and $L_{W}(-1)$ 
satisfying the same axioms for the corresponding operators coming from the vertex operator of the conformal 
element of a vertex operator algebra. We also refer the reader to Definition 3.1 in
\cite{H-twist-vo} for this notion in the special case that $V$ is a grading-restricted vertex algebra
 (not a superalgebra) and the automorphism of $V$ is $1_{V}$. Roughly speaking, a lower-bounded generalized $V$-module
is a $\C$-graded vector space $W=\coprod_{n\in \C}W_{[n]}$ equipped with a 
vertex operator map $Y_{W}: V\otimes W\to W[[x, x^{-1}]]$ and operators $L_{W}(0)$ and $L_{W}(-1)$
on $W$ satisfying all the 
axioms for an (ordinary) $V$-module except that for $n\in \C$, $W_{[n]}$ does not have to be finite dimensional and 
is the generalized eigenspace with the eigenvalue $n$ of $L_{W}(0)$ instead of the eigenspace with the  eigenvalue $n$
of $L_{W}(0)$. Module maps between lower-bounded generalized $V$-modules are defined 
in the obvious way as in Definition 1.1 
in \cite{H-cofiniteness}, not those defined in Definition 4.2 in \cite{H-exist-twisted-mod}. On the other hand, 
 if we replace $V$-module maps in the results below by those in Definition 4.2 in \cite{H-exist-twisted-mod},
these results still hold. The notion of generalized $V$-submodule of a lower-bounded generalized $V$-module 
is defined in the obvious way. A generalized $V$-submodule of a lower-bounded generalized $V$-module 
is certainly also lower bounded.

Let $W$ be a lower-bounded generalized $V$-module. 
For $n\in \N$, let 
$$\Omega_{n}(W)=\{w \in W \mid  (Y_{W})_{k}(v)w = 0\ {\rm for\ homogeneous}\ v\in V, 
\wt v - k - 1 < -n\}.$$
Then
$$\Omega_{n_{1}}(W)\subset \Omega_{n_{2}}(W)$$
for $n_{1}\le n_{2}$
and 
$$W=\bigcup_{n\in \N}\Omega_{n}(W).$$
So $\{\Omega_{n}(w)\}_{n\in \N}$ is an ascending filtration of $W$. 
Let 
$$Gr(W)=\sum_{n\in \N}Gr_{n}(W)$$
 be the 
associated graded space, where 
$$Gr_{n}(W)=\Omega_{n}(W)/\Omega_{n-1}(W).$$
Sometimes we shall use $[w]_{n}$ to denote the element $w+\Omega_{n-1}(W)$ of $Gr_{n}(W)$,
where $w\in \Omega_{n}(W)$. 

\begin{lemma}\label{Omega}
For $w\in \Omega_{n}(W)$ and $l\in \N$, 
$\res_{x}x^{l-1}Y_{W}(x^{L_{V}(0)}v, x)w\in \Omega_{n-l}(W)$.
\end{lemma}
\pf
The operator $\res_{x_{2}}x_{2}^{l-1}Y_{W}(x_{2}^{L_{V}(0)}v, x_{2})$ has weight $-l$.
Then for homogeneous $u\in V$, 
$(Y_{W})_{p}(u)\res_{x_{2}}x_{2}^{l-1}Y_{W}(x_{2}^{L_{V}(0)}v, x_{2})$ has weight 
$\wt u-p-1-l$. 
Consider the generalized $V$-submodule of $W$ generated by $w$. Then 
$(Y_{W})_{p}(u)\res_{x_{2}}x_{2}^{l-1}Y_{W}(x_{2}^{L_{V}(0)}v, x_{2})w$ is 
in this generalized $V$-submodule. Using the associativity for $Y_{W}$, we know that 
the generalized $V$-submodule generated by $w$ is 
spanned by elements of the form 
$(Y_{W})_{m}(\tilde{u})w$ for $\tilde{u}\in V$. So 
$(Y_{W})_{p}(u)\res_{x_{2}}x_{2}^{l-1}Y_{W}(x_{2}^{L_{V}(0)}v, x_{2})w$  
is a linear combination of such elements. But for homogeneous $w$, the weight of 
$(Y_{W})_{p}(u)\res_{x_{2}}x_{2}^{l-1}Y_{W}(x_{2}^{L_{V}(0)}v, x_{2})w$ 
is $\wt u-p-1-l+\wt w$. So the elements of the form 
$(Y_{W})_{m}(\tilde{u})w$ whose linear combination is 
$(Y_{W})_{p}(u)\res_{x_{2}}x_{2}^{l-1}Y_{W}(x_{2}^{L_{V}(0)}v, x_{2})w$ can also be chosen 
to be of weight $\wt u-p-1-l+\wt w$, that is, the weight $\wt \tilde{u}-m-1$ of $(Y_{W})_{m}(\tilde{u})$
is equal to $\wt u-p-1-l$. Since $w\in \Omega_{n}(W)$, $(Y_{W})_{m}(\tilde{u})w=0$ when 
$\wt \tilde{u}-m-1<-n$, or equivalently, $\wt u-p-1<-(n-l)$. So we have proved that 
$(Y_{W})_{p}(u)\res_{x_{2}}x_{2}^{l-1}Y_{W}(x_{2}^{L_{V}(0)}v, x_{2})w=0$  
when $\wt u-p-1<n-l$. This means that 
$\res_{x}x^{l-1}Y_{W}(x^{L_{V}(0)}v, x)w\in \Omega_{n-l}(W)$.
\epfv

By Lemma \ref{Omega}, the operator $\res_{x}x^{l-1}Y_{W}(x^{L_{V}(0)}v, x)$ in fact induces 
an operator, still denoted by the same notation, on $Gr(W)$, which maps $Gr_{n}(W)$ to 
$Gr_{n-l}(W)$. 

For $\mathfrak{v}=[v_{kl}]\in U^{\infty}(V)$, where $v_{kl}\in V$ and $k, l\in \N$, we define an operator
$\vartheta_{Gr(W)}(\mathfrak{v})$ on $Gr(W)$ as follows: For $\mathfrak{w}\in Gr(W)$,
we define
$$
\vartheta_{Gr(W)}(\mathfrak{v})\mathfrak{w}=\sum_{k, l\in \N}
\res_{x}x^{l-k-1}Y_{W}(x^{L_{V}(0)}v_{kl}, x)\pi_{Gr_{l}(W)}\mathfrak{w},
$$
where $\pi_{Gr_{l}(W)}$ is the projection from 
$Gr(W)$ to $Gr_{l}(W)$. 
Note that since $\mathfrak{w}$ is a sum of elements of $Gr_{l}(W)$ for finitely many $l\in \N$ and 
for each $l$, there are only finitely many nonzero $v_{kl}$, 
the sum over $k$ and $l$ is finite. So $\vartheta_{Gr(W)}(\mathfrak{v})\mathfrak{w}$ 
is indeed a well defined element of $Gr(W)$. 
In the case $\mathfrak{v}=[v]_{kl}$ and $\mathfrak{w}=[w]_{n}$
for $v\in V$, $w\in W$ and $k, l, n\in \N$, we have 
\begin{equation}\label{def-o-n-1}
\vartheta_{Gr(W)}([v]_{kl})[w]_{n}=\delta_{ln}
[\res_{x}x^{l-k-1}Y_{W}(x^{L_{V}(0)}v, x)w]_{k}.
\end{equation}
 In the case that $v$ is homogeneous and $w\in Gr_{l}(W)$, 
we have 
\begin{equation}\label{def-o-n-2}
\vartheta_{Gr(W)}([v]_{kl})[w]_{l}=[(Y_{W})_{\swt v+l-k-1}(v)w]_{k}.
\end{equation}

We now have a linear map 
\begin{align*}
\vartheta_{Gr(W)}: U^{\infty}(V)&\to \mbox{\rm End}\; Gr(W)\nn
\mathfrak{v}&\mapsto \vartheta_{Gr(W)}(\mathfrak{v}). 
\end{align*}
Let $Q^{\infty}(V)$ be the intersection of $\ker \vartheta_{Gr(W)}$ for all lower-bounded generalized $V$-modules $W$ and
$A^{\infty}(V)=U^{\infty}(V)/Q^{\infty}(V)$. 

We shall need the following lemma:

\begin{lemma}\label{L-1-der}
For $l\in \Z$, 
$k\in \N$ and $m\in \Z_{+}$  and $v\in V$, 
\begin{equation}\label{L-1-der-0}
\res_{x}x^{l-k-1}Y_{W}\left(x^{L_{V}(0)}\binom{L_{V}(-1)+L_{V}(0)+l}{k+m}v, x\right)=0.
\end{equation}
In particular, when $k=0$ and $m=1$, we have
\begin{equation}\label{L-1-der-0.5}
\res_{x}x^{l-1}Y_{W}\left(x^{L_{V}(0)}(L_{V}(-1)+L_{V}(0)+l)v, x\right)=0.
\end{equation}
For $l\in \Z$ and $v\in V$, 
\begin{equation}\label{L-1-der-0.7}
\res_{x}x^{l-k-1}Y_{W}\left(x^{L_{V}(0)}\binom{L_{V}(-1)+L_{V}(0)+l}{k}v, x\right)=
\res_{x}x^{l-k-1}Y_{W}(x^{L_{V}(0)}v, x).
\end{equation}
\end{lemma}
\pf 
For $l\in \Z$, $n\in \N$ and $v\in V$, 
using the $L(-1)$-derivative property for the vertex operator map $Y_{W}$ repeatedly, we have
\begin{equation}\label{L-1-der-1}
\frac{1}{n!}\frac{d^{n}}{dx^{n}}Y_{W}(x^{L_{V}(0)+l}v, x)=Y_{W}\left(x^{L_{V}(0)+l-n}
\binom{L_{V}(-1)+L_{V}(0)+l}{n}v, x\right).
\end{equation}
Multiplying $x^{p}$ to both sides and then taking $\res_{x}$, we obtain 
\begin{equation}\label{L-1-der-2}
\res_{x}\frac{x^{p}}{n!}\frac{d^{n}}{dx^{n}}Y_{W}(x^{L_{V}(0)+l}v, x)
=\res_{x}x^{l-n+p}Y_{W}\left(x^{L_{V}(0)}\binom{L_{V}(-1)+L_{V}(0)+l}{n}v, x\right).
\end{equation}
When $0\le p\le n-1$, the left-hand side of (\ref{L-1-der-2}) is $0$. Thus we obtain 
\begin{equation}\label{L-1-der-3}
\res_{x}x^{l-n+p}Y_{W}\left(x^{L_{V}(0)}\binom{L_{V}(-1)+L_{V}(0)+l}{n}v, x\right)=0.
\end{equation}
Let $n=k+m$ and $p=m-1$ in (\ref{L-1-der-3}) for 
$k\in \N$ and $m\in \Z_{+}$. Then we obtain (\ref{L-1-der-0}).

Let $p=-1$ and $n=k$ in (\ref{L-1-der-2}), we obtain 
\begin{equation}\label{L-1-der-4}
\res_{x}\frac{x^{-1}}{k!}\frac{d^{k}}{dx^{k}}Y_{W}(x^{L_{V}(0)+l}v, x)
=\res_{x}x^{l-k-1}Y_{W}\left(x^{L_{V}(0)}\binom{L_{V}(-1)+L_{V}(0)+l}{k}v, x\right).
\end{equation}
Since  left-hand side of (\ref{L-1-der-4}) is equal to 
$$\res_{x}x^{l-k-1}Y_{W}(x^{L_{V}(0)}v, x),$$
we obtain (\ref{L-1-der-0.7}).
\epfv

\begin{prop}\label{O-infty-mod}
We have  $O^{\infty}(V)\subset Q^{\infty}(V)$.
\end{prop}
\pf
We need to prove $\vartheta_{Gr(W)}(O^{\infty}(V))=0$ for every lower-bounded generalized $V$-module $W$.
For 
$$\res_{x_{0}}
x_{0}^{-k-l-p-2}(1+x_{0})^{l}[Y_{W}((1+x_{0})^{L(0)}v_{1}, x_{0})v_{2}]_{kl}
\in O^{\infty}(V),$$ 
where $v_{1}, v_{2}\in V$, $k, l, p\in \N$ and $w\in \Omega_{l}(W)$,  
we have
\begin{align}\label{O-infty-mod-1}
\vartheta_{Gr(W)}&(\res_{x_{0}}
x_{0}^{-k-l-p-2}(1+x_{0})^{l}[Y_{W}((1+x_{0})^{L(0)}v_{1},
x_{0})v_{2}]_{kl})[w]_{l}\nn
&=\res_{x_{2}}x_{2}^{l-k-1}\res_{x_{0}}
x_{0}^{-k-l-p-2}(1+x_{0})^{l}[Y_{W}(x_{2}^{L_{V}(0)}
Y_{V}((1+x_{0})^{L(0)}v_{1}, x_{0})v_{2}, x_{2})
\pi_{G_{l}(W)}w]_{k}\nn
&=\res_{x_{0}}\res_{x_{2}}
x_{0}^{-k-l-p-2}(1+x_{0})^{l}x_{2}^{l-k-1}\cdot\nn
&\quad\quad\quad\quad\quad\cdot 
[Y_{W}(Y_{V}(x_{2}^{L_{V}(0)}(1+x_{0})^{L(0)}v_{1}, x_{0}x_{2})
x_{2}^{L_{V}(0)}v_{2}, x_{2})w]_{k}\nn
&=\res_{x_{0}}\res_{x_{2}}
x_{0}^{-k-l-p-2}x_{2}^{-k-1}\res_{x_{1}}x_{1}^{l}x_{1}^{-1}
\delta\left(\frac{x_{2}+x_{0}x_{2}}{x_{1}}\right)\cdot\nn
&\quad\quad\quad\quad\quad\cdot 
[Y_{W}(Y_{V}(x_{1}^{L_{V}(0)}v_{1}, x_{0}x_{2})
x_{2}^{L_{V}(0)}v_{2}, x_{2})w]_{k}\nn
&=\res_{x_{0}}\res_{x_{2}}
x_{0}^{-k-l-p-2}x_{2}^{-k-1}\res_{x_{1}}x_{1}^{l}
x_{0}^{-1}x_{2}^{-1}
\delta\left(\frac{x_{1}-x_{2}}{x_{0}x_{2}}\right)\cdot\nn
&\quad\quad\quad\quad\quad\cdot 
[Y_{W}(x_{1}^{L(0)}v_{1}, x_{1})
Y_{W}(x_{2}^{L_{V}(0)}v_{2}, x_{2})w]_{k}\nn
&\quad -\res_{x_{0}}\res_{x_{2}}
x_{0}^{-k-l-p-2}x_{2}^{-k-1}\res_{x_{1}}x_{1}^{l}x_{0}^{-1}x_{2}^{-1}
\delta\left(\frac{x_{2}-x_{1}}{-x_{0}x_{2}}\right)\cdot\nn
&\quad\quad\quad\quad\quad\cdot 
[Y_{W}(x_{2}^{L_{V}(0)}v_{2}, x_{2})Y_{W}(x_{1}^{L(0)}v_{1}, x_{1})w]_{k}\nn
&=\res_{x_{1}}\res_{x_{2}}
x_{1}^{-k-p-2}(1-x_{1}^{-1}x_{2})^{-k-l-p-2}x_{2}^{l+p}
[Y_{W}(x_{1}^{L(0)}v_{1}, x_{1})
Y_{W}(x_{2}^{L_{V}(0)}v_{2}, x_{2})w]_{k}\nn
&\quad -\res_{x_{1}}\res_{x_{2}}
(-1+x_{1}x_{2}^{-1})^{-k-l-p-2}x_{1}^{l}x_{2}^{-k-2}
[Y_{W}(x_{2}^{L_{V}(0)}v_{2}, x_{2})Y_{W}(x_{1}^{L(0)}v_{1}, x_{1})w]_{k}.
\end{align}
Since $w\in \Omega_{l}(W)$ and the series $(1-x_{1}^{-1}x_{2})^{-k-l-p-2}$ contains only nonnegative powers of $x_{2}$,
$$\res_{x_{2}}(1-x_{1}^{-1}x_{2})^{-k-l-p-2}x_{2}^{l+p}Y_{W}(x_{2}^{L_{V}(0)}v_{2}, x_{2})w=0.$$
So the first term in the right-hand side of (\ref{O-infty-mod-1}) is $0$. 
Since $w\in \Omega_{l}(W)$ and the series $(-1+x_{1}x_{2}^{-1})^{-2k-p-2}$ contains only nonnegative powers of $x_{1}$, 
$$\res_{x_{1}}(-1+x_{1}x_{2}^{-1})^{-k-l-p-2}x_{1}^{l}Y_{W}(x_{1}^{L(0)}v_{1}, x_{1})
w=0.$$
So the second term in the right-hand side of (\ref{O-infty-mod-1}) is also $0$. 

Taking $l$ in (\ref{L-1-der-0.5}) to be $l-k$, we obtain 
\begin{align}
&\vartheta_{Gr(W)}([(L_{V}(-1)+L_{V}(0)+l-k)v]_{kl})[w]_{l}\nn
&\quad =[\res_{x}x^{l-k-1}Y_{W}\left(x^{L_{V}(0)}(L_{V}(-1)+L_{V}(0)+l-k)v, x\right)w]_{k}\nn
&\quad=0
\end{align}
for $v\in V$,  $k, l\in \N$ and $w\in \Omega_{l}(W)$. Thus we have $\vartheta_{Gr(W)}(O^{\infty}(V))=0$. 
\epfv

\begin{thm}\label{U-infty-mod}
Let $W$ be a lower-bounded generalized $V$-module. 
Then the linear map 
$$\vartheta_{Gr(W)}: U^{\infty}(V)\to \mbox{\rm End}\; Gr(W)$$
gives a
$U^{\infty}(V)$-module structure on $Gr(W)$ (that is, $\vartheta_{Gr(W)}$ 
is a homomorphism of (nonassociative) algebras from $U^{\infty}(V)$ to $\mbox{\rm End}\; Gr(W)$). 
In particular, $U^{\infty}(V)/\ker \vartheta_{Gr(W)}$
is an associative algebra isomorphic to a subalgebra of $\mbox{\rm End}\; Gr(W)$. 
\end{thm}
\pf
For $u, v\in V$, $k, n, l\in \N$ and $w\in \Omega_{l}(W)$, using (\ref{defn-diamond-1}), we have
\begin{align}\label{A-infty-mod-1}
&\vartheta_{Gr(W)}([u]_{kn}\diamond [v]_{nl})[w]_{l}\nn
&\quad=\res_{x_{0}}T_{k+l+1}((x_{0}+1)^{-k+n-l-1})(1+x_{0})^{l}\res_{x_{2}}
x_{2}^{l-k-1}\cdot\nn
&\quad\quad\quad\quad\quad\cdot 
[Y_{W}(x_{2}^{L_{V}(0)}Y_{V}((1+x_{0})^{L_{V}(0)}u, x_{0})v, x_{2})w]_{k}\nn
&\quad=\res_{x_{0}}\res_{x_{2}}T_{k+l+1}((x_{0}+1)^{-k+n-l-1})(1+x_{0})^{l}
x_{2}^{l-k-1}\cdot\nn
&\quad\quad\quad\quad\quad\cdot 
[Y_{W}(Y_{V}((x_{2}+x_{0}x_{2})^{L_{V}(0)}u, x_{0}x_{2})x_{2}^{L_{V}(0)}v, x_{2})w]_{k}\nn
&\quad=\res_{x_{0}}\res_{x_{2}}\res_{x_{1}}x_{1}^{-1}\delta\left(
\frac{x_{2}+x_{0}x_{2}}{x_{1}}\right)
T_{k+l+1}((x_{0}+1)^{-k+n-l-1})x_{1}^{l}
x_{2}^{-k-1}\cdot\nn
&\quad\quad\quad\quad\quad\cdot 
[Y_{W}(Y_{V}(x_{1}^{L_{V}(0)}u, x_{0}x_{2})x_{2}^{L_{V}(0)}v, x_{2})w]_{k}\nn
&\quad=\res_{x_{0}}\res_{x_{2}}\res_{x_{1}}x_{0}^{-1}x_{2}^{-1}
\delta\left(\frac{x_{1}-x_{2}}{x_{0}x_{2}}\right)T_{k+l+1}((x_{0}+1)^{-k+n-l-1})x_{1}^{l}x_{2}^{-k-1}\cdot\nn
&\quad\quad\quad\quad\quad\cdot 
[Y_{W}(x_{1}^{L_{V}(0)}u, x_{1})Y_{W}(x_{2}^{L_{V}(0)}v, x_{2})w]_{k}\nn
&\quad\quad -\res_{x_{0}}\res_{x_{2}}\res_{x_{1}}x_{0}^{-1}x_{2}^{-1}
\delta\left(\frac{x_{2}-x_{1}}{-x_{0}x_{2}}\right)T_{k+l+1}((x_{0}+1)^{-k+n-l-1})x_{1}^{l}x_{2}^{-k-1}\cdot\nn
&\quad\quad\quad\quad\quad\cdot 
[Y_{W}(x_{2}^{L_{V}(0)}v, x_{2})Y_{W}(x_{1}^{L_{V}(0)}u, x_{1})w]_{k}\nn
&\quad=\res_{x_{2}}\res_{x_{1}}T_{k+l+1}((x_{0}+1)^{-k+n-l-1})\lbar_{x_{0}=(x_{1}-x_{2})x_{2}^{-1}}
x_{1}^{l}x_{2}^{-k-2}\cdot\nn
&\quad\quad\quad\quad\quad\cdot 
[Y_{W}(x_{1}^{L_{V}(0)}u, x_{1})Y_{W}(x_{2}^{L_{V}(0)}v, x_{2})w]_{k}\nn
&\quad\quad -\res_{x_{2}}\res_{x_{1}}T_{k+l+1}((x_{0}+1)^{-k+n-l-1})
\lbar_{x_{0}=(-x_{2}+x_{1})x_{2}^{-1}}x_{1}^{l}x_{2}^{-k-2}\cdot\nn
&\quad\quad\quad\quad\quad\cdot 
[Y_{W}(x_{2}^{L_{V}(0)}v, x_{2})Y_{W}(x_{1}^{L_{V}(0)}u, x_{1})w]_{k}.
\end{align}
Since $w\in \Omega_{l}(W)$, the second term in the right-hand side of (\ref{A-infty-mod-1}) is $0$. 
Expanding $T_{k+l+1}((x_{0}+1)^{-k+n-l-1})$ explicitly, we see that the first term in the right-hand side of 
(\ref{A-infty-mod-1}) is equal to 
\begin{align}\label{A-infty-mod-2}
&\sum_{m=0}^{n}\binom{-k+n-l-1}{m}
\res_{x_{2}}\res_{x_{1}}(x_{1}-x_{2})^{-k+n-l-m-1}x_{2}^{k-n+l+m+1}
x_{1}^{l}x_{2}^{-k-2}\cdot\nn
&\quad\quad\quad\quad\quad\cdot 
[Y_{W}(x_{1}^{L_{V}(0)}u, x_{1})Y_{W}(x_{2}^{L_{V}(0)}v, x_{2})w]_{k}\nn
&\quad =\sum_{m=0}^{n}\sum_{j\in \N}\binom{-k+n-l-1}{m}\binom{-k+n-l-m-1}{j}(-1)^{j}
\res_{x_{2}}\res_{x_{1}}x_{1}^{-k+n-m-1-j}\cdot\nn
&\quad\quad\quad\quad\quad\cdot x_{2}^{-n+l+m-1+j}
[Y_{W}(x_{1}^{L_{V}(0)}u, x_{1})Y_{W}(x_{2}^{L_{V}(0)}v, x_{2})w]_{k}.
\end{align}
In the case $j>n-m$,
since $w\in \Omega_{l}(W)$, 
$$\wt v-(\wt v-n+l+m-1+j)-1<-l$$
when $v$ is homogeneous and hence we have 
$$\res_{x_{2}}x_{2}^{-n+l+m-1+j}Y_{W}(x_{2}^{L_{V}(0)}v, x_{2})w=0.$$
Hence those terms in the right-hand side of (\ref{A-infty-mod-2}) with $j>n-m$ is $0$.
So the right-hand side of (\ref{A-infty-mod-2}) is equal to 
\begin{align}\label{A-infty-mod-3}
&\sum_{m=0}^{n}\sum_{j=0}^{n-m}\binom{-k+n-l-1}{m}\binom{-k+n-l-m-1}{j}(-1)^{j}
\cdot\nn
&\quad\quad\quad\quad\quad\cdot  \res_{x_{2}}\res_{x_{1}}x_{1}^{-k+n-m-1-j}x_{2}^{-n+l+m-1+j}
[Y_{W}(x_{1}^{L_{V}(0)}u, x_{1})Y_{W}(x_{2}^{L_{V}(0)}v, x_{2})w]_{k}\nn
&\quad=\sum_{m=0}^{n}\sum_{p=m}^{n}\binom{-k+n-l-1}{m}\binom{-k+n-l-m-1}{p-m}(-1)^{p-m}
\cdot\nn
&\quad\quad\quad\quad\quad\cdot \res_{x_{2}}\res_{x_{1}}x_{1}^{-k+n-1-p} x_{2}^{-n+l-1+p}
[Y_{W}(x_{1}^{L_{V}(0)}u, x_{1})Y_{W}(x_{2}^{L_{V}(0)}v, x_{2})w]_{k}\nn
&\quad=\sum_{p=0}^{n}\left(\sum_{m=0}^{p}\binom{-k+n-l-1}{m}\binom{-k+n-l-m-1}{p-m}(-1)^{p-m}\right)
\cdot\nn
&\quad\quad\quad\quad\quad\cdot \res_{x_{2}}\res_{x_{1}}x_{1}^{-k+n-1-p} x_{2}^{-n+l-1+p}
[Y_{W}(x_{1}^{L_{V}(0)}u, x_{1})Y_{W}(x_{2}^{L_{V}(0)}v, x_{2})w]_{k}.
\end{align}
For $p=0, \dots, n$, 
\begin{align}\label{A-infty-mod-4}
&\sum_{m=0}^{p}\binom{-k+n-l-1}{m}\binom{-k+n-l-m-1}{p-m}(-1)^{p-m}\nn
&\quad=\sum_{m=0}^{p}\frac{(-k+n-l-1)\cdots (-k+n-l-m)}{m!}\cdot\nn
&\quad\quad\quad\quad\quad\cdot \frac{(-k+n-l-m-1)\cdots (-k+n-l-p)}{(p-m)!}(-1)^{p-m}\nn
&\quad=\sum_{m=0}^{p}\frac{(-k+n-l-1)\cdots  (-k+n-l-p)}{p!}\frac{p!}{m!(p-m)!}(-1)^{p-m}\nn
&\quad=\binom{-k+n-l-1}{p}\sum_{m=0}^{p}\binom{p}{m}(-1)^{p-m}\nn
&\quad=\binom{-k+n-l-1}{p}(-1+1)^{p}\nn
&\quad =\binom{-k+n-l-1}{p}\delta_{p, 0}.
\end{align}
Using (\ref{A-infty-mod-4}), we see that the right-hand side of (\ref{A-infty-mod-3}) is equal to 
\begin{align}\label{A-infty-mod-5}
&\res_{x_{2}}\res_{x_{1}}x_{1}^{-k+n-1} x_{2}^{-n+l-1}
[Y_{W}(x_{1}^{L_{V}(0)}u, x_{1})Y_{W}(x_{2}^{L_{V}(0)}v, x_{2})w]_{k}\nn
&\quad =\vartheta_{Gr(W)}([u]_{kn})
[\res_{x_{2}}x_{2}^{l-n-1}Y_{W}(x_{2}^{L_{V}(0)}v, x_{2})w]_{n}\nn
&\quad =\vartheta_{Gr(W)}([u]_{kn})\vartheta_{Gr(W)}([v]_{nl})[w]_{l}.
\end{align}
From (\ref{A-infty-mod-1}), (\ref{A-infty-mod-2}), (\ref{A-infty-mod-3}) and (\ref{A-infty-mod-5}),
we obtain 
$$\vartheta_{Gr(W)}([u]_{kn}\diamond [v]_{nl})=\vartheta_{Gr(W)}([u]_{kn})\vartheta_{Gr(W)}([v]_{nl})$$
for $u, v\in V$ and $k, n, l\in \N$. Thus $\vartheta_{Gr(W)}$ gives an $U^{\infty}(V)$-module structure on $W$. 
\epfv

\begin{lemma}
Let $L_{U}(-1)$ and $L_{U}(0)$ be operators on a vector space $U$ satisfying 
$$[L_{U}(0), L_{U}(-1)]=L_{U}(-1).$$
We have 
\begin{equation}\label{L-1-L0}
e^{xL_{U}(-1)}(1+x)^{L_{U}(0)}=(1+x)^{L_{U}(-1)+L_{U}(0)}.
\end{equation}
\end{lemma}
\pf 
This can be proved easily by showing
$$\frac{d}{dx}e^{xL_{U}(-1)}(1+x)^{L_{U}(0)}(1+x)^{-(L_{U}(-1)+L_{U}(0))}=0$$
so that it must be independent of $x$ and then setting $x=0$ to obtain 
$$e^{xL_{U}(-1)}(1+x)^{L_{U}(0)}(1+x)^{-(L_{U}(-1)+L_{U}(0))}=1_{U}.$$
\epfv

We can now write down explicitly the expressions of 
 elements of the form $[v]_{kl}\diamond \one^{\infty}$ for $v\in V$ and 
$k, l\in \N$ satisfying $k\le l$.

\begin{lemma}
For $v\in V$ and $k, l\in \N$,
\begin{equation}\label{one-on-right}
[v]_{kl}\diamond\one^{\infty}= \sum_{m=0}^{l}\binom{-k-1}{m}\left[\binom{L_{V}(-1)+L_{V}(0)+l}{k+m}v\right]_{kl}.
\end{equation}
\end{lemma}
\pf
By the definition (\ref{defn-diamond}) of $\diamond$ and the skew-symmetry of $Y_{V}$, 
\begin{align}\label{one-on-right-1}
([v]_{kl}\diamond_{N}\one^{\infty})_{mn}&=\delta_{km}\delta_{ln}\res_{x}T_{k+l+1}((x+1)^{-k+n-l-1})
(1+x)^{n}Y_{V}((1+x)^{L(0)}v, x)\one\nn
&=\delta_{km}\delta_{ln}\res_{x}T_{k+l+1}((x+1)^{-k+n-l-1})
(1+x)^{n}e^{xL_{V}(-1)}(1+x)^{L_{V}(0)}v.
\end{align}
Thus we obtain 
\begin{equation}\label{one-on-right-2}
[v]_{kl}\diamond\one^{\infty}=\res_{x}T_{k+l+1}((x+1)^{-k-1})
(1+x)^{l}[e^{xL_{V}(-1)}(1+x)^{L_{V}(0)}v]_{kl}.
\end{equation}
Using (\ref{L-1-L0}) with $U=V$, expanding the formal series explicitly and then evaluating the 
formal residue, we see that the right-hand side of (\ref{one-on-right-2}) is equal to 
\begin{align}\label{one-on-right-3}
&\res_{x}T_{k+l+1}((x+1)^{-k-1})
(1+x)^{l}[(1+x)^{L_{V}(0)+L_{V}(0)}v]_{kl}\nn
&\quad =\sum_{m=0}^{l}\binom{-k-1}{m}\res_{x}
x^{-k-m-1}[(1+x)^{L_{V}(-1)+L_{V}(0)+l}v]_{kl}\nn
&\quad =\sum_{m=0}^{l}\sum_{j\in \N}\binom{-k-1}{m}\res_{x}x^{-k-m-1+j}\binom{L_{V}(-1)+L_{V}(0)+l}{j}v\nn
&\quad =\sum_{m=0}^{l}\binom{-k-1}{m}\binom{L_{V}(-1)+L_{V}(0)+k}{k+m}v,
\end{align}
proving  (\ref{one-on-right}).
\epfv

\begin{prop}\label{right-id}
For $v\in V$ and $k, l\in \N$, 
$[v]_{kl}\diamond\one^{\infty}-[v]_{kl}\in O^{\infty}(V)$. For $u, v\in V$ and $k, l, n\in \N$, 
$([v]_{kl}\diamond\one^{\infty}-[v]_{kl})\diamond [u]_{ln}\in O^{\infty}(V)$.
\end{prop}
\pf
For $m\in \N$, 
\begin{align}\label{right-id-1}
\binom{L_{V}(-1)+L_{V}(0)+l}{k+m}v&=\binom{(L_{V}(-1)+L_{V}(0)+l-k)+k}{k+m}v\nn
&= \binom{k}{k+m}v+(L_{V}(-1)+L_{V}(0)+l-k)\tilde{v}_{m}\nn
&\equiv \left\{
\begin{array}{ll}0&m\in \Z_{+}\\
v&m=0
\end{array}\right. \mod O^{\infty}(V),\nn
\end{align}
where $\tilde{v}_{m}$ is an element of $V$ depending on $m$. 
Thus by (\ref{one-on-right}), 
$$[v]_{kl}\diamond\one^{\infty}\equiv v\mod O^{\infty}(V).$$

By (\ref{one-on-right}), (\ref{right-id-1}) and (\ref{der-O-infty}), 
\begin{align*}
([v]_{kl}\diamond\one^{\infty})\diamond [u]_{ln}
&=\sum_{m=0}^{l}\binom{-k-1}{m}\left[\binom{L_{V}(-1)+L_{V}(0)+l}{k+m}v\right]_{kl}\diamond
 [u]_{ln}\nn
&=\sum_{m=0}^{l}\binom{-k-1}{m} \binom{k}{k+m}[v]_{kl}\diamond
 [u]_{ln}\nn
& \quad +\sum_{m=0}^{l}\binom{-k-1}{m}\left[(L_{V}(-1)+L_{V}(0)+l-k)\tilde{v}_{m}\right]_{kl}\diamond
 [u]_{ln}\nn
&\equiv [v]_{kl}\diamond
 [u]_{ln} \mod O^{\infty}(V).
\end{align*}
\epfv

\begin{thm}\label{assoc-alg}
The product $\diamond$ on $U^{\infty}(V)$ induces a product, denoted still by $\diamond$, on 
$A^{\infty}(V)=U^{\infty}(V)/Q^{\infty}(V)$ such that $A^{\infty}(V)$  equipped with $\diamond$ is an associative algebra
with $\one^{\infty}+Q^{\infty}(V)$ as identity. 
Moreover, the associated graded space $Gr(W)$ of the ascendant
filtration $\{\Omega_{n}(W)\}_{n\in \N}$ of a lower-bounded generalized $V$-module $W$  is an 
$A^{\infty}(V)$-module. 
\end{thm}
\pf
Since $\ker \vartheta_{Gr(W)}$ for a lower-bounded generalized $V$-module 
$W$ is a two-sided ideal of $U^{\infty}(V)$, $Q^{\infty}(V)$ as the intersection of 
such two-sided ideals is still a two-sided ideal of $U^{\infty}(V)$. Thus 
$\diamond$ on $U^{\infty}(V)$ induces a product on
$A^{\infty}(V)$. Since for each lower-bounded generalized $V$-module $W$,
the quotient algebra $U^{\infty}(V)/\ker \vartheta_{Gr(W)}$ is associative, we have 
$$\mathfrak{v}_{1}\diamond (\mathfrak{v}_{2}\diamond \mathfrak{v}_{3})-
\mathfrak{v}_{1}\diamond (\mathfrak{v}_{2}\diamond \mathfrak{v}_{3})\in \ker \vartheta_{Gr(W)}$$ 
for 
$\mathfrak{v}_{1}, \mathfrak{v}_{2},  \mathfrak{v}_{3}\in U^{\infty}(V)$.
Then we have 
$$\mathfrak{v}_{1}\diamond (\mathfrak{v}_{2}\diamond \mathfrak{v}_{3})-
\mathfrak{v}_{1}\diamond (\mathfrak{v}_{2}\diamond \mathfrak{v}_{3})\in \bigcap_{W} \ker \vartheta_{Gr(W)}
=Q^{\infty}(V)$$  
for 
$\mathfrak{v}_{1}, \mathfrak{v}_{2},  \mathfrak{v}_{3}\in U^{\infty}(V)$.
Thus
$A^{\infty}(V)=U^{\infty}(V)/Q^{\infty}(V)$
is also associative. 

By definition, we have $\one^{\infty}\diamond [v]_{kl}=[v]_{kl}$.
So $\one^{\infty}$ is in fact a left identity of the algebra $U^{\infty}(V)$.
By Proposition \ref {O-infty-mod}, $O^{\infty}(V)\subset Q^{\infty}(V)$. Then by Proposition \ref{right-id},
we have 
$$([v]_{kl}+Q^{\infty}(V))\diamond (\one^{\infty}+Q^{\infty}(V))=[v]_{kl}+Q^{\infty}(V).$$
So $\one^{\infty}+Q^{\infty}(V)$ is an identity of $A^{\infty}(V)$. 

For a lower-bounded generalized $V$-module $W$, by Theorem \ref{U-infty-mod},
$Gr(W)$ is a module for $U^{\infty}(V)/\ker \vartheta_{Gr(W)}$. 
Since $Q^{\infty}(V)$ is a two-sided subideal of $\ker \vartheta_{Gr(W)}$, $Gr(W)$ is an 
$A^{\infty}(V)$-module. 
\epfv

The ideal $Q^{\infty}(V)$ of $U^{\infty}(V)$ is defined using all lower-bounded generalized $V$-modules.
From Prposition 3.3 in \cite{H-aa-int-op}, it is now known that besides elements of
$O^{\infty}(V)$, $Q^{\infty}(V)$ also contains elements corresponding to the Jacobi identity 
for $V$. We conjecture that  $Q^{\infty}(V)$ is generated by $O^{\infty}(V)$ and these elements. 

The only result above on $O^{\infty}(V)$ and $Q^{\infty}(V)$  is Proposition \ref{O-infty-mod}. 
Below we give another result on $O^{\infty}(V)$ (Proposition \ref{l-r-prod-L-1L0}). 
To prove this result, we need the following commutator formula:

\begin{lemma}
For $v \in V$,
\begin{equation}\label{comm-L-1+L-0}
[L_{V}(-1)+L_{V}(0), Y_{V}((1+x)^{L_{V}(0)}v, x)]
=Y_{V}((1+x)^{L_{V}(0)}(L_{V}(-1)+L_{V}(0))v, x).
\end{equation}
\end{lemma}
\pf
By the $L(-1)$ and $L(0)$-commutator formula with the vertex operator map $Y_{V}$
and the fact that the weight of $L_{V}(-1)$ is $1$, 
\begin{align*}
&[L_{V}(-1)+L_{V}(0), Y_{V}((1+x)^{L_{V}(0)}v, x)]\nn
&\quad=Y_{V}(((1+x)L_{V}(-1)+L_{V}(0))(1+x)^{L_{V}(0)}v, x)\nn
&\quad=Y_{V}((1+x)^{L_{V}(0)}(L_{V}(-1)+L_{V}(0))v, x).
\end{align*}

By Theorem \ref{U-infty-mod}, every lower-bounded generalized $V$-module is an 
$A^{\infty}(V)$-module. 
\epfv

\begin{prop}\label{l-r-prod-L-1L0}
For $u, v\in V$ and $k, n, l\in V$, 
both 
$$[(L_{V}(-1)+L_{V}(0)+n-k)u]_{kn}\diamond [v]_{nl}$$
 and 
$$[v]_{kn}\diamond[(L_{V}(-1)+L_{V}(0)+l-n)u]_{nl}$$ 
are in $O^{\infty}(V)$.
\end{prop}
\pf
For $u, v\in V$, $k, l, m\in \N$, 
by definition, 
\begin{align}\label{der-O-infty-1}
&[(L_{V}(-1)+L_{V}(0)+n-k)u]_{kn}\diamond [v]_{nl}\nn
&\quad=\res_{x}T_{k+l+1}((x+1)^{-k+n-l-1})(1+x)^{l}[Y_{V}((1+x)^{L_{V}(0)}
(L_{V}(-1)+L_{V}(0)+n-k)u, x)v]_{kl}\nn
&\quad=\res_{x}T_{k+l+1}((x+1)^{-k+n-l-1})(1+x)^{k-n+l+1}\frac{d}{dx}[Y_{V}((1+x)^{L_{V}(0)+n-k}u, x)v]_{kl}\nn
&\quad=-\res_{x}\left(\frac{d}{dx}T_{k+l+1}((x+1)^{-k+n-l-1})(1+x)^{k-n+l+1}\right)
[Y_{V}((1+x)^{L_{V}(0)+n-k}u, x)v]_{kl}\nn
&\quad=-\res_{x}\Biggl(\left(\frac{d}{dx}T_{k+l+1}((x+1)^{-k+n-l-1})\right)(1+x)^{k-n+l+1}\nn
&\quad\quad\quad\quad\quad\quad\quad\quad\quad\quad
+(k-n+l+1)T_{k+l+1}((x+1)^{-k+n-l-1})(1+x)^{k-n+l}\Biggr)\cdot\nn
&\quad\quad\quad\quad\quad \cdot (1+x)^{-k+n}[Y_{V}((1+x)^{L_{V}(0)}u, x)v]_{kl}.
\end{align}
Applying $-\frac{1+x}{k-n+l+1}\frac{d}{dx}$ to bother sides of (\ref{exp-x+1-k-1}),
we obtain 
\begin{align}\label{der-O-infty-2}
&(x+1)^{-k+n-l-1}\nn
&\; =-\frac{1+x}{k-n+l+1}\frac{d}{dx}T_{k+l+1}((x+1)^{-k+n-l-1})
-\frac{1+x}{k-n+l+1}\frac{d}{dx}R_{k+l+1}((x+1)^{-k+n-l-1}).
\end{align}
Since the first and second terms in the right-hand side of (\ref{der-O-infty-2}) contain 
only the terms with powers in $x^{-1}$ less than or equal to and larger than, respectively, $k+l+2$,
we must have 
$$-\frac{1+x}{k-n+l+1}\frac{d}{dx}T_{k+l+1}((x+1)^{-k+n-l-1})=T_{k+l+2}((x+1)^{-k+n-l-1}),$$
or equivalently, 
\begin{align}\label{der-O-infty-3}
&(1+x)\frac{d}{dx}T_{k+l+1}((x+1)^{-k-1})\nn
&\quad =-(k-n+l+1)T_{k+l+1}((x+1)^{-k+n-l-1})-(k-n+l+1)\binom{-k+n-l-1}{n+1}x^{-k-l-2},
\end{align}
Using (\ref{der-O-infty-3}), the right-hand side of (\ref{der-O-infty-1}) becomes
$$(k-n+l+1)\binom{-k+n-l-1}{n+1}\res_{x}x^{-k-l-2}(1+x)^{l}[Y_{V}((1+x)^{L_{V}(0)}u, x)v]_{kl} \in O^{\infty}(V).$$
Thus we obtain 
\begin{align}\label{der-O-infty}
&[(L_{V}(-1)+L_{V}(0)+n-k)u]_{kn}\diamond [v]_{nl}\nn
&\quad =(k-n+l+1)\binom{-k+n-l-1}{n+1}\res_{x}x^{-k-l-2}(1+x)^{l}[Y_{V}((1+x)^{L_{V}(0)}u, x)v]_{kl} \nn
&\quad \in O^{\infty}(V).
\end{align}

For $u, v\in V$, $k, l, n\in \N$ satisfying $n\le k\le l$, 
by the definition,  (\ref{comm-L-1+L-0}) and $l-n=(l-k)-(n-k)$
\begin{align}\label{assoc-alg-1}
&[v]_{kn}\diamond[(L_{V}(-1)+L_{V}(0)+l-n)u]_{nl}\nn
&\quad =\res_{x}T_{k+l+1}((x+1)^{-k+n-l-1})(1+x)^{l}\cdot\nn
&\quad\quad\quad\quad\cdot 
[Y_{V}((1+x)^{L_{V}(0)}v, x)(L_{V}(-1)+L_{V}(0)+l-n)u]_{kl}\nn
&\quad =\res_{x}T_{k+l+1}((x+1)^{-k+n-l-1})(1+x)^{l}\cdot\nn
&\quad\quad\quad\quad\cdot 
[(L_{V}(-1)+L_{V}(0)+l-k)Y_{V}((1+x)^{L_{V}(0)}v, x)u]_{kl}\nn
&\quad \quad -\res_{x}T_{k+l+1}((x+1)^{-k+n-l-1})(1+x)^{l}\cdot\nn
&\quad\quad\quad\quad\cdot 
[Y_{V}((1+x)^{L_{V}(0)}(L_{V}(-1)+L_{V}(0)+n-k)v, x)u]_{kl}.
\end{align}
The first term in the right-hand side of (\ref{assoc-alg-1}) is by definition in 
$O^{\infty}(V)$. The second term in the right-hand side of (\ref{assoc-alg-1})
is equal to $[(L_{V}(-1)+L_{V}(0)+n-k)v]_{kn}\diamond [u]_{nl}$, which is 
also in $O^{\infty}(V)$ (\ref{der-O-infty}). So 
$$[v]_{kn}\diamond[(L_{V}(-1)+L_{V}(0)+l-n)u]_{nl}\in O^{\infty}(V).$$
\epfv

\begin{rema}
{\rm In \cite{DJ}, for $m, n, p\in\Z_{+}$, a product 
$*_{m, p}^{n}$ on $V$ is introduced. In terms of these products, we have
$[u]_{kn}\diamond[v]_{nl}=[u*_{l, n}^{k}v]_{kl}$. 
For each $m, n\in \Z_{+}$, a subspace $O'_{m, n}(V)$ is also introduced in \cite{DJ}. 
In terms of these subspaces, $O^{\infty}(V)$ can be easily shown to be spanned 
by infinite linear combinations of 
elements of the form $[v]_{kl}$ for $v\in O'_{k, l}(V)$ with each pair $(k, l)$ appearing in the 
linear combinations only finitely many times. For each $m, n\in \Z_{+}$, 
a subspace $O_{m, n}(V)$ of $V$ containing in particular $O'_{m, n}(V)$ and associators of the products
$*_{r, p}^{q}$ is further introduced. By taking a suitable subspace of the direct product of 
the quotient spaces $A_{m, n}(V)=V/O_{m, n}(V)$ for $m, n\in \N$, it is possible to 
obtain an associative algebra. 
This associative algebra can be identified as the quotient of the nonassociative algebra
$U^{\infty}(V)$ by the ideal generated by $O^{\infty}(V)$ and all the associators of the product $\diamond$.
Such a construction of associative algebras works for any nonassociative algebra with a subspace. 
As we have mentioned in the introduction, the algebra 
$A^{\infty}(V)$ plays the role of universal enveloping algebra of $V$ for
the category of lower-bounded generalized  $V$-modules. 
So we now use an analogy with Lie algebra to compare the associative algebra that one can get from 
\cite{DJ} and the associative algebra $A^{\infty}(V)$ introduced in this paper. 
Given a Lie algebra, one can generate a free nonassociative algebra.
Then the quotient of this free associative algebra by the ideal generated by
associators is an associative algebra isomorphic to the tensor algebra 
generated by the Lie algebra. To obtain the universal envelpoping algebra of the Lie algebra, 
we need to further take the quotient by the Jacobi identity and the skew-symmetry relations. 
The nonassociative algebra $U^{\infty}(V)$ is analogous to the free nonassociative algebra generated by the 
Lie algebra. The associative algebra that one can get from \cite{DJ} is analogous to the tensor algebra of the Lie algebra. 
The associative algebra $A^{\infty}(V)$ is analogous to the universal enveloping algebra of the Lie algebra. 
Certainly $A^{\infty}(V)$ is a quotient of the associative algebra that one can get from \cite{DJ}.
But $A^{\infty}(V)$ is not equal to this associative algebra. 
In fact, Corollary 3.6 in \cite{H-aa-int-op} says that $Q^{\infty}(V)$ contains the elements 
\begin{align*}
&\sum_{\stackrel{j\in \N}{n+p-j\ge 0}}(-1)^{j}\binom{p}{j}[v]_{k,n+p-j}\diamond
[u]_{n+p-j,l+p}\nn
&\quad -\sum_{\stackrel{j\in \N}{l-n+k+p-j\ge 0}}(-1)^{p-j}\binom{p}{j}
[u]_{k,l-n+k+p-j}\diamond [v]_{l-n+k+p-j,l+p}\nn
&\quad-\sum_{j\in \N}\binom{\wt v+n-k-1}{j}
[(Y_{V})_{p+j}(v)u]_{k,l+p}
\end{align*}
for $k, l, n\in \N$, $p\in \Z$ such that $l+p\in \N$, $u\in V$ and homogeneous $v\in V$,
corresponding to elements giving the Jacobi identity. Such elements 
are in general not in $O_{k, l+p}(V)$ in \cite{DJ}.  We conjecture that $Q^{\infty}(V)$ is generated by $O^{\infty}(V)$ and 
these elements.}
\end{rema}

\setcounter{equation}{0} \setcounter{thm}{0} 
\section{Lower-bounded generalized $V$-modules and graded $A^{\infty}$-modules}

We study the relations between lower-bounded generalized $V$-modules 
and suitable $A^{\infty}(V)$-modules in this section. 

Note that $W$ is graded by the generalized eigenspaces of $L_{W}(0)$. Since $\Omega_{n}(W)$ for $n\in \N$
is invariant under $L_{W}(0)$, $L_{W}(0)$ induces an operator on $Gr_{n}(W)=\Omega_{n}(W)/\Omega_{n-1}(W)$ 
such that $Gr_{n}(W)$ is also graded by  the generalized eigenspaces of this operator. 
These operators on $Gr_{n}(W)$ for $n\in \N$ together define an operator, denoted by $L_{Gr(W)}(0)$, on $Gr(W)$ 
preserving the $\N$-grading on $Gr(W)$. 
Then $Gr(W)$ is also graded by the generalized eigenspaces of $L_{Gr(W)}(0)$. 

For $v\in V$, $k\in \Z$ and $w\in \Omega_{n}(W)$, by the $L(-1)$-commutator formula, 
$$(Y_{W})_{k}(v)L_{W}(-1)w=L_{W}(-1)(Y_{W})_{k}(v)w+k(Y_{W})_{k-1}(v)w.$$
When $\wt v-k-1<-(n+1)$, we have $\wt v-k-1<-n$ and $\wt v-(k-1)-1<-n$. 
So in this case, $(Y_{W})_{k}(v)w=(Y_{W})_{k-1}(v)w=0$ since $w\in \Omega_{n}(W)$.
Thus $(Y_{W})_{k}(v)L_{W}(-1)w=0$ when $\wt v-k-1<-(n+1)$. This means that 
$L_{W}(-1)w\in \Omega_{n+1}(W)$. In particular, $L_{W}(-1)$ induces a linear map 
from $Gr_{n}(W)$ to $Gr_{n+1}(W)$ for $n\in \N$. These maps for $n\in \N$ together
define an operator, denoted by $L_{Gr(W)}(-1)$, on $Gr(W)$. 

The operators $L_{Gr(W)}(0)$,
$L_{Gr(W)}(-1)$ and $\vartheta_{Gr(W)}([v]_{kl})$ 
satisfy the same commutator formulas as those between $L_{W}(0)$, $L_{W}(-1)$ and 
$\res_{x}x^{l-k-1}Y_{W}(x^{L_{V}(0)}v, x)$ for $v\in V$ and $k, l\in \N$.
These structures on $Gr(W)$ motivates the following definition: 

\begin{defn}\label{N-graded-A-inf-mod}
{\rm Let $G$ be an $A^{\infty}(V)$-module with the $A^{\infty}(V)$-module structure on $G$ given by a homomorphism
$\vartheta_{G}: A^{\infty}(V)\to {\rm End}\; G$ of 
associative algebras. We say that $G$ is a {\it graded $A^{\infty}(V)$-module}
if the following conditions are 
satisfied:

\begin{enumerate}

\item $G$ is graded by $\N$, that is, $G=\coprod_{n\in \N}G_{n}$, and for $v\in V$, $k, l\in \N$,
$\vartheta_{G}([v]_{kl}+Q^{\infty}(V))$ maps $G_{n}$ to 
$0$ when $n\ne l$ and to $G_{k}$ when $n=l$. 

\item $G$ is a direct sum of generalized eigenspaces of  an operator $L_{G}(0)$ on $G$, $G_{n}$ for $n\in \N$ 
are invariant under $L_{G}(0)$ and
the real parts of the eigenvalues of $L_{G}(0)$ have a lower bound. 

\item There is an operator $L_{G}(-1)$ on $G$ mapping $G_{n}$ to $G_{n+1}$ for $n\in \N$.

\item The commutator relations
\begin{align*}
{[L_{G}(0), L_{G}(-1)]}&=L_{G}(-1),\nn
{[L_{G}(0), \vartheta_{G}([v]_{kl}+Q^{\infty}(V))]}&=(k-l)\vartheta_{G}([v]_{kl}+Q^{\infty}(V)),\nn
{[L_{G}(-1), \vartheta_{G}([v]_{kl}+Q^{\infty}(V))]}&=\vartheta_{G}([L_{V}(-1)v]_{(k+1)l}+Q^{\infty}(V))
\end{align*}
hold for $v\in V$ and 
$k, l\in \N$

\end{enumerate}
A graded $A^{\infty}(V)$-algebra $G$ is said to be {\it nondegenerate} if it satisfies in addition the following condition:
For $\textsl{g}\in G_{l}$, if $\vartheta_{G}([v]_{0l}+Q^{\infty}(V))\textsl{g}=0$ for all $v\in V$, then $\textsl{g}=0$. 
Let $G_{1}$ and $G_{2}$ be graded $A^{\infty}(V)$-modules.
A {\it graded $A_{\infty}(V)$-module map} from 
$G_{1}$ to $G_{2}$ is an $A^{N}(V)$-module map $f: G_{1}\to G_{2}$ such that 
$f((G_{1})_{n})\subset (G_{2})_{n}$,  
$f\circ L_{G_{1}}(0)=L_{G_{2}}(0)\circ f$ and $f\circ L_{G_{1}}(-1)=L_{G_{2}}(-1)\circ f$.
A {\it graded $A^{\infty}(V)$-submodule} of a graded $A^{\infty}(V)$-module $G$ is an $A^{\infty}(V)$-submodule 
 of $G$ that is also an $\N$-graded subspace of $G$ and invariant under 
the operators $L_{G}(0)$ and $L_{G}(-1)$.
A graded $A^{\infty}(V)$-module $G$ is said to be {\it generated by a subset $S$} if 
$G$  is equal  to the smallest graded $A^{\infty}(V)$-submodule containing $S$, or equivalently, 
$G$  is spanned by homogeneous elements with respect to the $\N$-grading and the grading given by 
$L_{G}(0)$ obtained by applying elements of $A^{\infty}(V)$, $L_{G}(0)$ and $L_{G}(-1)$ to 
homogeneous summands of elements of $S$. 
 A graded $A^{\infty}(V)$-module is said to be 
{\it irreducible} if it has no nonzero proper graded $A^{\infty}(V)$-submodules. A graded $A^{\infty}(V)$-module 
is said to be {\it completely reducible} if it is a direct sum of irreducible graded $A^{\infty}(V)$-modules.}
\end{defn}

From Theorem \ref{assoc-alg} and 
the properties of a lower-bounded generalized $V$-module $W$ and its associated graded 
space $Gr(W)$, we obtain immediately:

\begin{thm}\label{W-adm}
For a lower-bounded generalized $V$-module $W$, $Gr(W)$ is a nondegenerate graded $A^{\infty}(V)$-module.
Let $W_{1}$ and $W_{2}$ be lower-bounded generalized $V$-modules 
and $f: W_{1}\to W_{2}$ a $V$-module map. Then $f$ induces a graded $A^{\infty}(V)$-module map
$Gr(f): Gr(W_{1})\to Gr(W_{2})$. 
\end{thm}

We now give a direct and explicit description of $Gr(W)$ for a completely 
reducible lower-bounded generalized $V$-module $W$. 
In this case, 
$$W=\coprod_{\mu\in \mathcal{M}}W^{\mu},$$
where $\mathcal{M}$ is an index set and $W^{\mu}$ for $\mu \in \mathcal{M}$ 
are irreducible lower-bounded generalized $V$-modules. For $\mu\in \mathcal{M}$,
since $W^{\mu}$ is irreducible, there exists $h^{\mu}\in \C$ such that
$$W^{\mu}=\coprod_{n\in \N}W^{\mu}_{[h^{\mu}+n]},$$
where as usual, $W^{\mu}_{[h^{\mu}+n]}$ for $n\in \N$ is the subspace of $W^{\mu}$ 
of weight $h^{\mu}+n$, and  $W^{\mu}_{[h^{\mu}]}\ne 0$.
For $n\in \N$, let
$$G_{n}(W)=\coprod_{\mu\in \mathcal{M}}W_{[h^{\mu}+n]}.$$
Then 
$$W=\coprod_{n\in \N}G_{n}(W).$$
For $n\in \N$, let 
$$T_{n}(W)=\coprod_{m=0}^{n}G_{m}(W).$$
It is clear that $T_{n}(W)\subset \Omega_{n}(W)$. In particular, $G_{n}(W)\subset \Omega_{N}(W)$ for $n\le N$.
Let $e_{W}: W\to Gr(W)$ be defined by $e_{W}(w)=w+\Omega_{n-1}(W)$ for $w\in G_{n}(W)$ and $n\in \N$. 
Then $e_{W}$ preserves the $\N$-grading. 
We also define a map $\vartheta_{W}: U^{\infty}(V)\to {\rm End}\;W$ by 
$$\vartheta_{W}(\mathfrak{v})w=\sum_{k, l\in \N}
\res_{x}x^{l-k-1}Y_{W}(x^{L_{V}(0)}v_{kl}, x)\pi_{G_{l}(W)}w
$$
for $\mathfrak{v}\in U^{\infty}$ and $w\in W$, 
where $\pi_{G_{l}(W)}$ is the projection from 
$W$ to $G_{l}(W)$. 
In the case $\mathfrak{v}=[v]_{kl}$ and $w\in G_{n}(W)$
for $v\in V$ and $k, l, n\in \N$, we have 
\begin{equation}\label{W-def-o-n-1}
\vartheta_{W}([v]_{kl})w=\delta_{ln}
\res_{x}x^{l-k-1}Y_{W}(x^{L_{V}(0)}v, x)w.
\end{equation}

\begin{prop}\label{W=GrW}
Let $W$ be a completely reducible lower-bounded generalized $V$-module. Then
$\Omega_{n}(W)=T_{n}(W)$ for $n\in \N$. Moreover, 
$W$ equipped with $\vartheta_{W}$
 is a nondegenerate graded $A^{\infty}(V)$-module and 
$e_{W}: W\to Gr(W)$ is an  isomorphism of graded $A^{\infty}(V)$-modules. 
\end{prop}
\pf
If $T_{n}(W)\ne \Omega_{n}(W)$, then there exists homogeneous $w\in \Omega_{n}(W)$ but $w\not\in T_{n}(W)$. 
Then $w=\sum_{\mu\in \mathcal{M}}w^{\mu}$, where 
$w_{\mu}\in W^{\mu}$ for $\mu\in \mathcal{M}$ and only finitely many $w^{\mu}$ is not $0$. 
Since $w$ is homogeneous, we can assume that $w^{\mu}$ for $\mu \in \mathcal{M}$ are  homogeneous.
Since $w\in \Omega_{n}(W)$ but $w\not\in T_{n}(W)$, there is at least 
one $w^{\mu}$ such that $w^{\mu}\in \Omega_{n}(W^{\mu})$ but $w^{\mu}\not\in T_{n}(W^{\mu})
=\coprod_{m=0}^{n}W^{\mu}_{[h^{\mu}+m]}$.
Let $W_{0}^{\mu}$ be the 
generalized $V$-submodule of $W^{\mu}$
generated by such a $w^{\mu}$.  Since $w^{\mu}\not\in T_{n}(W^{\mu})$, $w^{\mu}\ne 0$ and hence
$W_{0}^{\mu}\ne 0$. But $W^{\mu}$ is irreducible. So $W_{0}^{\mu}=W^{\mu}$. 
Since $w^{\mu}$ is homogeneous, there is $m\in \N$ such that $\wt w^{\mu}\in W^{\mu}_{[h^{\mu}+m]}$.
Since $w^{\mu}\not\in T_{n}(W^{\mu})$, we must have $m>n$. Since $W^{\mu}=W_{0}^{\mu}$,
$W^{\mu}$ is spanned by elements of the form $(Y_{W})_{k}(v)w^{\mu}$ for $v\in V$ and $k\in \Z$.
Since $w^{\mu}\in \Omega_{n}(W^{\mu})$, $(Y_{W})_{k}(v)w^{\mu}=0$ for 
homogeneous $v\in V$ and $k\in \Z$ satisfying $\wt v-k-1<-n$. Thus the homogeneous subspaces of 
$W^{\mu}_{[h^{\mu}+m-n-p]}=0$ for $p\in \Z_{+}$. 
But for $p=m-n\in \Z_{+}$, $W^{\mu}_{[h^{\mu}+m-n-p]}=W^{\mu}_{[h^{\mu}]}\ne 0$. 
Contradiction. Thus $T_{n}(W)=\Omega_{n}(W)$. 

For $n\in \N$, we have $Gr_{n}(W)=\Omega_{n}(W)/\Omega_{n-1}(W)=T_{n}(W)/T_{n-1}(W)$.
Then $e_{W}\mbar_{G_{n}(W)}$ is clearly a linear isomorphism from $G_{n}(W)$ to $T_{n}(W)/T_{n-1}(W)=Gr_{n}(W)$. 
This shows that $e_{W}$ is an isomorphism of graded spaces. 
For $v\in V$, $k, l\in \N$ and $w\in G_{l}(W)$, 
\begin{align*}
e_{W}(\vartheta_{W}([v]_{kl})w)&=e_{W}(\res_{x}x^{l-k-1}Y_{W}(x^{L_{V}(0)}v, x)w)\nn
&=\res_{x}x^{l-k-1}Y_{W}(x^{L_{V}(0)}v, x)w+T_{k-1}(W)\nn
&=\vartheta_{Gr(W)}([v]_{kl})e_{W}(w).
\end{align*}
Thus we have $e_{W}\circ \vartheta_{W}=\vartheta_{Gr(W)}\circ e_{W}$. 
In particular, the $A^{\infty}(V)$-module structure on $Gr(W)$ given by $\vartheta_{Gr(W)}$ 
is transported to $W$ by $e_{W}$ so that $W$ equipped with $\vartheta_{W}$
 is an $A^{\infty}(V)$-module and 
$e_{W}: W\to Gr(W)$ is an  isomorphism of $A^{\infty}(V)$-modules. 
\epfv

\begin{thm}\label{reduct-A-infty}
A lower-bounded generalized $V$-module $W$ is irreducible or completely reducible
 if and only if the nondegenerate graded
$A^{\infty}(V)$-module $Gr(W)$  is irreducible or completely reducible, respectively. 
\end{thm}
\pf
Let $W$ be an irreducible lower-bounded generalized $V$-module. By Theorem  \ref{W=GrW},
$W$ is a nondegenerate graded $A^{\infty}(V)$-module isomorphic to $Gr(W)$. Let $W_{0}$ be a nonzero graded $A^{\infty}(V)$-submodule 
of the graded $A^{\infty}(V)$-module $W$. For a homogeneous element $v\in V$, $n\in \Z$
and $w\in W_{0}$, 
$$\res_{x}x^{n}Y_{W}(v, x)w=\sum_{l\in \N}\vartheta_{W}([v]_{(\swt v-n-1+l)l})\pi_{G_{l}(W)}w\in W_{0}.$$ 
This means that $W_{0}$ is invariant under the action of the vertex operators on $W$. By the definition of 
graded  $A^{\infty}(V)$-submodule, $W_{0}$ is invariant under the actions of $L_{W}(0)$ and $L_{W}(-1)$ and
is the direct sum of generalized eigenspaces of $L_{W}(0)\mbar_{W_{0}}$.
Thus $W_{0}$ is also a nonzero lower-bounded generalized $V$-submodule of $W$. 
Since $W$ is an  irreducible lower-bounded generalized $V$-module, $W_{0}=W$. 
So as a graded $A^{\infty}(V)$-module, $W$ is also irreducible. Since 
as a graded $A^{\infty}(V)$-module, $Gr(W)$ is equivalent to $W$, we see that 
$Gr(W)$ is irreducible. 

Conversely, 
assume that for a lower-bounded generalized  $V$-module $W$, the nondegenerate 
graded $A^{\infty}(V)$-module $Gr(W)$ is 
irreducible. Let $W_{0}$ be  a nonzero generalized $V$-submodule of $W$. Then $\Omega_{n-1}(W_{0})\subset 
\Omega_{n-1}(W)$ for $n\in \N$ (when $n=0$, $\Omega_{-1}(W)=0$).  We have a map from $Gr(W_{0})$ to $Gr(W)$ given by 
$w_{0}+\Omega_{n-1}(W_{0}) \mapsto w_{0}+\Omega_{n-1}(W)$ for $n\in \N$ and $w_{0}\in \Omega_{n}(W_{0})$. 
This map is an injective graded $A^{\infty}(V)$-module map. So the image of $Gr(W_{0})$ under this map 
is a graded $A^{\infty}(V)$-submodule of $Gr(W)$. Since $W_{0}$ is nonzero, $Gr(W_{0})$ is nonzero.
Since $Gr(W)$ is irreducible,  the image of $Gr(W_{0})$ under this map is equal to 
$Gr(W)$. Now it is easy to derive $W_{0}=W$. 
In fact, for $n\in \N$,  the image of $Gr_{n}(W_{0})$ under the map
from $Gr(W_{0})$ to $Gr(W)$ above   is 
$\{w_{0}+\Omega_{n}(W)\mid w_{0}\in \Omega_{n}(W_{0})\}$. So 
$Gr_{n}(W)=\{w_{0}+\Omega_{n-1}(W)\mid w_{0}\in \Omega_{n}(W_{0})\}$.
For $n=0$, we obtain $\Omega_{0}(W)=Gr_{0}(W)=Gr(W_{0})=\Omega_{0}(W_{0})$. Assume that 
$\Omega_{n-1}(W)=\Omega_{n-1}(W_{0})$. Given $w\in \Omega_{n}(W)$,
$w+\Omega_{n-1}(W)\in Gr_{n}(W)$. By
$Gr_{n}(W)=\{w_{0}+\Omega_{n-1}(W)\mid w_{0}\in \Omega_{n}(W_{0})\}$,
there exists $w_{0}\in \Omega_{n}(W_{0})$ such that $w+\Omega_{n-1}(W)=w_{0}+\Omega_{n-1}(W)$,
or equivalently, $w-w_{0}\in \Omega_{n-1}(W)=\Omega_{n-1}(W_{0})$. 
Thus $w\in \Omega_{n}(W_{0})$. This shows $\Omega_{n}(W)=\Omega_{n}(W_{0})$ for $n\in \N$.
Then we have $W=\cup_{n\in \N}\Omega_{n}(W)=\cup_{n\in \N}\Omega_{n}(W_{0})=W_{0}$. 
So $W$ is irreducible. 

Assume that a lower-bounded generalized $V$-module $W$ is completely reducible. 
Then $W=\coprod_{\mu\in \mathcal{M}}W^{\mu}$,
where $W^{\mu}$ for $\mu\in \mathcal{\mu}$ are irreducible generalized $V$-modules. From what we have proved above, 
$W^{\mu}$ for $\mu\in \mathcal{M}$ as graded $A^{\infty}(V)$-modules are also irreducible. 
So $W$ as a graded $A^{\infty}(V)$-module is completely reducible.
But $Gr(W)$ is equivalent to $W$ as a graded $A^{\infty}(V)$-module by Proposition \ref{W=GrW}. So
$Gr(W)$  is also completely reducible. Conversely, 
assume that for a lower-bounded generalized $V$-module $W$, the graded $A^{\infty}(V)$-module 
$Gr(W)$ is completely reducible. Then $Gr(W)=\coprod_{\mu\in \mathcal{M}}G^{\mu}$,
where $G^{\mu}$ for $\mu\in \mathcal{M}$ are irreducible nondegenerate 
graded $A^{\infty}(V)$-submodules of $Gr(W)$.
For $\mu\in \mathcal{M}$, since $G^{\mu}$ is  a nondegenerate graded $A^{\infty}(V)$-submodule of $Gr(W)$,
we have $G^{\mu}_{n}\subset 
Gr_{n}(W)=\Omega_{n}(W)/\Omega_{n-1}(W)$. Let $W^{\mu}$ be the subspace of $W$ consisting of 
elements of the form $w^{\mu}\in \Omega_{n}(W)$ such that $w^{\mu}+\Omega_{n-1}(W)\in G^{\mu}_{n}$
for $n\in \N$. Since $G^{\mu}$ is a nondegenerate graded $A^{\infty}(V)$-submodule of $Gr(W)$, 
for $v\in V$, $k, l\in \N$ and $w^{\mu}\in \Omega_{l}(W)$ such that $w^{\mu}+\Omega_{l-1}(W)\in G^{\mu}_{l}$,
$$\res_{x}x^{l-k-1}Y_{W}(x^{L_{V}(0)}v, x)w^{\mu}+\Omega_{k-1}(W)\in G^{\mu}_{k}.$$
By the definition of $W^{\mu}$, we obtain 
$\res_{x}x^{l-k-1}Y_{W}(x^{L_{V}(0)}v, x)w^{\mu}\in W^{\mu}$. 
Since $w^{\mu}\in \Omega_{l}(W)$, $\res_{x}x^{l-k-1}Y_{W}(x^{L_{V}(0)}v, x)w^{\mu}=0$
for $k\in -\Z_{+}$. Thus 
$\res_{x}x^{l-k-1}Y_{W}(x^{L_{V}(0)}v, x)w^{\mu}\in W^{\mu}$ for $k\in \N$ are all the nonzero 
coefficients of $Y_{W}(v, x)w^{\mu}$. This means that $W^{\mu}$ is closed under the action of 
the vertex operators on $W$. Since $G^{\mu}$ is invariant under the actions of $L_{Gr(W)}(0)$ and $L_{Gr(W)}(-1)$ 
and is a direct sum of generalized eigenspaces of $L_{Gr(W)}(0)$,
$W^{\mu}$ is invariant under the actions of $L_{W}(0)$ and $L_{W}(-1)$ 
and is a direct sum of generalized eigenspaces of $L_{W}(0)$. Thus $W^{\mu}$ is a generalized $V$-submodule of $W$. 

Let $w^{\mu}+\Omega_{n-1}(W^{\mu})\in Gr_{n}(W^{\mu})$,  where $n\in \N$ and $w^{\mu}\in \Omega_{n}(W^{\mu})
\subset \Omega_{n}(W)$. By the definition of $W^{\mu}$, we see that since $w^{\mu}$ is
 an element of $W^{\mu}$, 
$w^{\mu}+\Omega_{n-1}(W)\in G^{\mu}_{n}$. So we obtain a linear map from 
$Gr_{n}(W^{\mu})$ to $G^{\mu}_{n}$ given by $w^{\mu}+\Omega_{n-1}(W^{\mu})\mapsto
w^{\mu}+\Omega_{n-1}(W)$ for $w^{\mu}+\Omega_{n-1}(W^{\mu})\in Gr_{n}(W^{\mu})$. 
These maps for $n\in \N$ give a map from $Gr(W^{\mu})$ to $G^{\mu}$. 
It is clear that this map is a graded $A^{\infty}(V)$-module map. 
If the image $w^{\mu}+\Omega_{n-1}(W)$ of $w^{\mu}+\Omega_{n-1}(W^{\mu})\in Gr_{n}(W^{\mu})$
under this map is $0$ in $G^{\mu}$, then $w^{\mu}\in \Omega_{n-1}(W)$.
But $w^{\mu}\in \Omega_{n}(W^{\mu})\subset W^{\mu}$. So $w^{\mu}\in \Omega_{n-1}(W^{\mu})$
and $w^{\mu}+\Omega_{n-1}(W^{\mu})$ is $0$ in $Gr(W^{\mu})$. This means that this 
graded $A^{\infty}(V)$-module map is injective. In particular, the image of $Gr(W^{\mu})$ under this map 
is a nonzero nondegenerate graded $A^{\infty}(V)$-submodule of $G^{\mu}$. But $G^{\mu}$ is irreducible.
So $Gr(W^{\mu})$ must be equivalent to $G^{\mu}$ and is therefore also irreducible. 
From what we have proved above, since $Gr(W^{\mu})$ is irreducible, 
$W^{\mu}$ is irreducible. This shows that $W$ is complete reducible. 
\epfv

Theorem \ref{reduct-A-infty} implies that there is a map from the set of the equivalence classes of 
 irreducible lower-bounded generalized $V$-modules to the set of equivalence classes of 
 irreducible nondegenerate graded $A^{\infty}(V)$-modules. 
This map is in fact a bijection. To prove this, we need to construct a lower-bounded generalized $V$-module $S(G)$ 
from a nondegenerate graded $A^{\infty}(V)$-module $G$. 
We use the construction in Section 5 of \cite{H-const-twisted-mod}. 
Take the generating fields for the grading-restricted vertex algebra $V$ to be 
$Y_{V}(v, x)$ for $v\in V$. By definition, 
$G$ is a direct sum of generalized eigenspaces of $L_{G}(0)$ and the real parts of the eigenvalues of 
$L_{G}(0)$ has a lower bound $B\in \R$. We take $M$ and $B$ in Section 5 of \cite{H-const-twisted-mod}
to be $G$ and the lower bound $B$ above. 
Using the construction in  Section 5 of \cite{H-const-twisted-mod},
we obtain a universal lower-bounded generalized $V$-module 
$\widehat{G}^{[1_{V}]}_{B}$. For simplicity, we shall denote it simply by 
$\widehat{G}$. 

By Theorem 3.3 in \cite{H-exist-twisted-mod} and the construction in  Section 5 of \cite{H-const-twisted-mod}
and by identifying elements of the form 
$(\psi_{\widehat{G}}^{a})_{-1, 0}$ with basis elements $\g^{a}\in G$ for $a\in A$
for a basis $\{\g^{a}\}_{a\in A}$ of $G$, 
we see that $\widehat{G}$ is generated by $G$ (in the sense of Definition 3.1 in \cite{H-exist-twisted-mod}). 
Moreover,  after identifying $(\psi_{\widehat{G}}^{a})_{-1, 0}$ with basis elements $w^{a}\in G$ for $a\in A$,
Theorems 3.3 and 3.4 in \cite{H-exist-twisted-mod} in fact
say that elements of the form 
$L_{\widehat{G}}(-1)^{p}w^{a}$ for $p\in \N$ and $a\in A$ are linearly independent and 
$\widehat{G}$ is spanned by elements obtained by applying the components of the vertex operators
to these elements.  In particular, $G$ can be embedded into $\widehat{G}$ as 
a subspace. So from now on, we shall view $G$ as a subspace of $\widehat{G}$.
Let $J_{G}$  be  the generalized $V$-submodule 
of $\widehat{G}$ generated by elements of the forms
\begin{equation}\label{J-G-0}
\res_{x}x^{l-k-1}Y_{\widehat{G}}(x^{L_{V}(0)}v, x)\g
\end{equation}
for $ l\in \N$, $k\in -\Z_{+}$ and $\g\in G_{l}$, 
\begin{equation}\label{J-G}
\res_{x}x^{l-k-1}Y_{\widehat{M}}(x^{L_{V}(0)}v, x)\g-\vartheta_{G}([v]_{kl}+Q^{\infty}(V))\g
\end{equation}
for $v\in V$, $k, l\in \N$, $\g\in G_{l}$ and 
\begin{equation}\label{J-G-2}
L_{\widehat{G}}(-1)\g-L_{G}(-1)\g
\end{equation}
for $l\in \N$, $\g\in G_{l}$. 

Let $S(G)=\widehat{G}/J_{G}$. Then $S(G)$ is a lower-bounded generalized $V$-module. 
Let $\pi_{S(G)}$ be the projection from $\widehat{G}$ to $S(G)$. 
Since $\widehat{G}$ is generated by $G$ (in the sense of Definition 3.1 in \cite{H-exist-twisted-mod}),
$S(G)$ is generated by $\pi_{S(G)}(G)$ (in the same sense). In particular, 
$S(G)$ is spanned by elements of the form 
\begin{equation}\label{S-G}
\res_{x}x^{(l+p)-n-1}Y_{S(G)}(x^{L_{S(G)}(0)}v, x)L_{S(G)}(-1)^{p}\pi_{S(G)}(\g)
\end{equation}
for $v\in V$, $n, l, p\in \N$ and $\g\in G_{l}$. For $n\in \N$, let $G_{n}(S(G))$ be the subspace of 
$S(G)$ spanned by elements of the form (\ref{S-G}) for $v\in V$, $l, p\in \N$ and $\g\in G_{l}$.

\begin{prop}\label{S-G-prop}
Let $G$ be a  graded $A^{\infty}(V)$-module. 

\begin{enumerate}

\item For $n\in \N$, $G_{n}(S(G))=\pi_{S(G)}(G_{n})$ and for $n_{1}\ne n_{2}$, $G_{n_{1}}(S(G))\cap G_{n_{2}}(S(G))=0$. 
Moreover, ${\displaystyle S(G)=\pi_{S(G)}(G)=\coprod_{n\in \N}G_{n}(S(G))}$.

\item If $G$ is nondegenerate, then for $n\in \N$, ${\displaystyle  \Omega_{n}(S(G))=
\coprod_{j=0}^{n}\pi_{S(G)}(G_{j})=\coprod_{j=0}^{n}G_{j}(S(G))}$.

\item If $G$ is nondegenerate, then $Gr(S(G))$ is equivalent to $G$ as a graded $A^{\infty}(V)$-module.

\end{enumerate}
\end{prop}
\pf
Since elements of the forms (\ref{J-G}) and (\ref{J-G-2}) are in $J_{G}$, for $n\in \N$, 
the element (\ref{S-G}) for $v\in V$ for $l, p\in \N$ and $\g\in G_{l}$ is in fact equal to 
\begin{equation}\label{S-g-2}
\pi_{S(G)}(\vartheta_{G}([v]_{n(l+p)}+Q^{\infty}(V))L_{G}(-1)^{p}\g).
\end{equation}
Since $\vartheta_{G}([v]_{n(l+p)}+Q^{\infty}(V))L_{G}(-1)^{p}\g$ for $l, p\in \N$ and $\g\in G_{l}$  
certainly span $G_{n}$ and elements of the form (\ref{S-G})  for $v\in V$ for $l, p\in \N$ and $\g\in G_{l}$
span $G_{n}(S(G))$, elements of the form (\ref{S-g-2}) for $v\in V$ for $l, p\in \N$ and $\g\in G_{l}$ also span
$G_{n}(S(G))$. Thus $G_{n}(S(G))=\pi_{S(G)}(G_{n})$. 
When $n_{1}\ne n_{2}$, we know $G_{n_{1}}\cap G_{n_{2}}=0$.
Then $G_{n_{1}}(S(G))\cap G_{n_{2}}(S(G))=\pi_{S(G)}(G_{n_{1}}\cap G_{n_{2}})=0$.
As is mentioned above,
$S(G)$ is spanned by elements of the form (\ref{S-G})
for $v\in V$, $k, l, p\in \N$ and $\g\in G_{l}$. But we already see that 
(\ref{S-G}) is in fact equal to (\ref{S-g-2}). Thus $S(G)=\pi_{S(G)}(G)$.
Since $G_{n}(S(G))=\pi_{S(G)}(G_{n})$ and $G_{n_{1}}(S(G))\cap G_{n_{2}}(S(G))=0$, 
we have $S(G)=\pi_{S(G)}(G)=\coprod_{n\in \N}G_{n}(S(G))$.

By definition, for $j\le n$, $G_{j}(S(G))\subset \Omega_{n}(S(G))$. 
Then for $j=0, \dots, n$, $\pi_{S(G)}(G_{j})\subset  \Omega_{j}(S(G))\subset \Omega_{n}(S(G))$. 
So we obtain $\pi_{S(G)}(\coprod_{j=0}^{n}G_{j})\subset  \Omega_{n}(S(G))$.
If $G$ is nondegenerate, nonzero  elements of 
$G_{j}$ for $j>n$ are not in $\Omega_{n}(\widehat{G})$. 
From the construction of $\widehat{G}$, nonzero elements of the form (\ref{J-G-0}),  (\ref{J-G})
or  (\ref{J-G-2})  are not in $G\subset \widehat{G}$.
In particular,  the intersection of $J(G)$ with $G$ is $0$.
So $\pi_{S(G)}\mbar_{G}: G\to S(G)$ is injective.
Since $\pi_{S(G)}\mbar_{G}$ is injective, 
we conclude that nonzero  elements of 
$ \pi_{S(G)}(G_{j})$ for $j>n$ are not in $\Omega_{n}(S(G))$. 
So we have  
$$\Omega_{n}(S(G))=\pi_{S(G)}\left(\coprod_{j=0}^{n}G_{j}\right)=
\coprod_{j=0}^{n}\pi_{S(G)}(G_{j})=\coprod_{j=0}^{n}G_{j}(S(G)).$$

Since $\Omega_{n}(S(G))=\coprod_{j=0}^{n}G_{j}(S(G))$ for $n\in \N$ when $G$ is nondegenerate, 
we see that  as a $\N$-graded space, $Gr(S(G))$ is isomorphic to 
$\coprod_{n\in \N}G_{n}(S(G))=\pi_{S(G)}(G)$.
We use $f_{G}$ to denote the isomorphism from $Gr(S(G))$ to $\pi_{S(G)}(G)$. Then we have 
$$f_{G}\circ \vartheta_{Gr(S(G))}([v]_{kl}+Q^{\infty}(V))
=\res_{x}x^{l-k-1}Y_{S(G)}(x^{L_{S(G)}(0)}v, x)\circ f_{G}$$
for $v\in V$, $k, l\in \N$,
$f_{G}\circ L_{Gr(S(G))}(0)=L_{S(G)}(0)\circ 
f_{G}$ and $f_{G}\circ L_{Gr(S(G))}(-1)=L_{S(G)}(-1)\circ f_{G}$.
We have proved that $\pi_{S(G)}\mbar_{G}$ is injective and surjective and preserves the $\N$-gradings.
So it is an isomorphism of graded spaces from $G$ to $S(G)$. From the fact that $\pi_{S(G)}$ is 
a $V$-module map and on $G\subset \widehat{G}$,
$L_{\widehat{G}}(0)=L_{G}(0)$ and $L_{\widehat{G}}(-1)=L_{G}(-1)$, we have
$$\pi_{S(G)}\mbar_{G}\circ \vartheta_{G}([v]_{kl}+Q^{\infty}(V))
=\res_{x}x^{l-k-1}Y_{S(G)}(x^{L_{S(G)}(0)}v, x)\circ \pi_{S(G)}\mbar_{G}$$
for $v\in V$, $k, l\in \N$,
$\pi_{S(G)}\mbar_{G}\circ L_{G}(0)=L_{S(G)}(0)\circ 
\pi_{S(G)}\mbar_{G}$ and $\pi_{S(G)}\mbar_{G}\circ L_{G}(-1)=L_{S(G)}(-1)\circ \pi_{S(G)}\mbar_{G}$.
Then by the properties of $f_{G}$ and $\pi_{S(G)}\mbar_{G}$ above, 
we see that $(\pi_{S(G)}\mbar_{G})^{-1}\circ f_{G}$
is an equivalence of graded $A^{\infty}(V)$-modules from 
$Gr(S(G))$ to $G$.
\epfv

\begin{rema}
{\rm Note that our construction of the lower-bounded generalized $V$-module 
$\widehat{G}$ seems to depend on the lower bound $B$ of the real parts of the 
eigenvalues of $L_{G}(0)$. But by Proposition \ref{S-G-prop}, $S(G)$ depends only on 
$G$, not on $B$.}
\end{rema}

\begin{thm}\label{irred-A-inf-m-bij}
The set of the equivalence classes of irreducible lower-bounded generalized $V$-modules is in bijection with 
the set of  the equivalence classes of irreducible nondegenerate graded $A^{\infty}(V)$-modules. 
\end{thm}
\pf
Let $[\mathfrak{W}]_{\rm irr}$ be the set  of the equivalence classes of irreducible lower-bounded generalized 
$V$-modules and $[\mathfrak{G}]_{\rm irr}$  the 
set of  the equivalence classes of irreducible nondegenerate graded $A^{\infty}(V)$-modules.
Given an irreducible lower-bounded generalized $V$-module $W$, by Theorem \ref{reduct-A-infty}, 
$Gr(W)$ is an irreducible nondegenerate graded $A^{N}(V)$-module. Thus we obtain a map $f:
[\mathfrak{W}]_{\rm irr}\to [\mathfrak{G}]_{\rm irr}$ 
given by $f([W])=[Gr(W)]$, where $[W]\in [\mathfrak{W}]_{\rm irr}$ is the equivalence class containing
the irreducible  lower-bounded generalized $V$-module $W$ and $[Gr(W)]\in [\mathfrak{G}]_{\rm irr}$ is the equivalence 
class containing the irreducible nondegenerate graded $A^{N}(V)$-module $Gr(W)$. By Proposition \ref{W=GrW},
$[Gr(W)]=[W]$ in $[\mathfrak{G}]_{\rm irr}$, where $W$ is viewed as a nondegenerate graded $A^{\infty}(V)$-module. 

Given an irreducible nondegenerate graded $A^{\infty}(V)$-module $G$, we have a lower-bounded generalized 
$V$-module $S(G)$. By Proposition \ref{S-G-prop}, $Gr(S(G))$ is equivalent to $G$. Since $G$ is irreduible, 
$Gr(S(G))$ is also irreducible. Then by Theorem \ref{reduct-A-infty},  $S(G)$ is an
irreducible lower-bounded generalized  $V$-module. 
Thus we obtain a map 
$g: [\mathfrak{G}]_{\rm irr}\to [\mathfrak{W}]_{\rm irr}$ given by $g([G])=[S(G)]$. 

We still need to show that $f$ and $g$ are inverse to each other. 
By Proposition \ref{S-G-prop},
$Gr(S(G))$ is equivalent to $G$ for an irreducible nondegenerate graded $A^{\infty}(V)$-module $G$. 
We obtain $[Gr(S(G))]=[G]$. This means 
$f(g([G]))=[G]$. So we have $f\circ g=1_{[\mathfrak{G}]_{\rm irr}}$. 

Let $W$ be an irreducible  lower-bounded generalized $V$-module. By Theorem \ref{reduct-A-infty},
$Gr(W)$ is an irreducible nondegenerate 
graded $A^{\infty}(V)$-module. We then have a lower-bounded generalized  $V$-module 
$S(Gr(W))$.  By Proposition \ref{S-G-prop}, $Gr(S(Gr(W)))$ is equivalent to 
$Gr(W)$ as a graded $A^{\infty}(V)$-module. Since $Gr(W)$ is irreducible, $Gr(S(Gr(W)))$
is also irreducible. By Theorem \ref{reduct-A-infty}, $S(Gr(W))$ is an irreducible lower-bounded generalized  $V$-module. 
Since both $W$ and $S(Gr(W))$ are irreducible,  by  Proposition \ref{W=GrW}, $W$ and $S(Gr(W))$
are nondegenerate 
graded $A^{\infty}(V)$-modules and are equivalent to $Gr(W)$ and $Gr(S(Gr(W)))$, respectively. 
But we already know that $Gr(S(Gr(W)))$ is equivalent to 
$Gr(W)$ as a graded $A^{\infty}(V)$-module. So both $W$ and $S(Gr(W))$ are 
equivalent to $Gr(W)$ as graded $A^{\infty}(V)$-modules. Since vertex operators on 
$W$ and $S(Gr(W))$ can be expressed using the 
actions of elements of $A^{\infty}(V)$, we see that $W$ and $S(Gr(W))$ are also equivalent as 
lower-bounded generalized $V$-modules. Thus $[S(Gr(W))]=[W]$, 
or $g(f([W]))=[W]$. So  $g\circ f=1_{[\mathfrak{W}]_{\rm irr}}$. 
\epfv

\setcounter{equation}{0} \setcounter{thm}{0} 
\section{Subalgebras of $A^{\infty}(V)$}

We give some very special subalgebras of $A^{\infty}(V)$ and prove that they are isomorphic to 
the Zhu algebra $A(V)$ \cite{Z}  and its generalizations $A_{N}(V)$ for 
$N\in \N$ by Dong-Li-Mason 
\cite{DLM} in Subsection 4.1. 
Then we introduce the main interesting and new subalgebras $A^{N}(V)$ for $N\in \N$  of 
$A^{\infty}(V)$ in Subsection 4.2. Note that we use the superscript $N$ instead of the subscript $N$
to distinguish this algebra from $A_{N}(V)$ in \cite{DLM}. 

\subsection{Zhu algebra and the generalizations by Dong-Li-Mason}

Let $$U_{00}(V)=\{[v]_{00}\;|\;v\in V\}
\subset U^{\infty}(V).$$
Then $U_{00}(V)$ can be canonically identified with 
$V$ through the map $i_{00}: U_{00}(V)\to V$ given by $i_{00}([v]_{00})=v$ for $v\in V$.
Since  by  (\ref{defn-diamond-1}), 
$$[u]_{00}\diamond [v]_{00}
=\res_{x}x^{-1}\left[Y_{V}((1+x)^{L(0)}u, x)v\right]_{00},$$
$U_{00}(V)$ is closed under the product $\diamond$. Let 
$$A_{00}(V)=\{[v]_{00}+Q^{\infty}(V)\;|\;v\in V\}.$$

\begin{thm}\label{zhu-isom}
The subspace $A_{00}(V)$ of $A^{\infty}(V)$ is closed under $\diamond$ and is thus a subalgebra
of $A^{\infty}(V)$ with 
$[\one]_{00}+Q^{\infty}(V)$ as its identity. 
The associative algebra $A_{00}(V)$ 
 is isomorphic to the Zhu algebra $A(V)$ in \cite{Z} and, in particular, $[\omega]_{00}+
Q^{\infty}(V)$ is in the center of $A_{00}(V)$ if $V$ is a vertex operator algebra with 
the conformal vector $\omega$. 
\end{thm}

Since this result is a special case of the result on the generalizations $A_{N}(V)$ 
in \cite{DLM}, we will not give a proof. The proof is the special case $N=0$ of the proof of 
Theorem \ref{DLM-isom} below for $A_{N}(V)$.

Let 
$$U^{NN}=\left\{\sum_{k=0}^{N}[v]_{kk}\;\lbar
\;v\in V\right\}
\subset U^{\infty}(V).$$
By  the definition of $\diamond$, 
\begin{align*}
\left(\sum_{k=0}^{N}[u]_{kk}\right)\diamond 
\left(\sum_{k=0}^{N}[v]_{kk}\right)
& =\sum_{k, l=0}^{N}[u]_{kk}\diamond [v]_{ll}\nn
&  =\sum_{k=0}^{N}[u]_{kk}\diamond [v]_{kk}
\nn
& =\sum_{k=0}^{N}
\res_{x}T_{2k+1}((x+1)^{-k-1})(1+x)^{k}\left[Y_{V}((1+x)^{L(0)}u, x)v\right]_{kk}.
\end{align*}
Let
$$A^{NN}(V)=\left\{\sum_{k=0}^{N}[v]_{kk}+Q^{\infty}(V)\;\lbar
\;v\in V\right\}
\subset A^{\infty}(V).$$
Note that $A^{00}(V)=A_{00}(V)$. 
Also let 
\begin{align*}
\one^{N}&=\sum_{k=0}^{N}[\one]_{kk},\\
\omega^{N}&=\sum_{k=0}^{N}[\omega]_{kk},
\end{align*}

\begin{thm}\label{DLM-isom}
The subspace $A^{NN}(V)$ of $A^{\infty}(V)$ is closed under $\diamond$ and is thus a subalgebra
of $A^{\infty}(V)$ with $\one^{N}+Q^{\infty}(V)$ as the identity.
The associative algebra $A^{NN}(V)$
 is isomorphic to the associative algebra $A_{N}(V)$ of Dong, Li and Mason  in \cite{DLM} and,
in particular,  
$\omega^{N}+Q^{\infty}(V)$ is in the center of $A^{NN}(V)$ if $V$ is a vertex operator algebra with 
the conformal vector $\omega$.
\end{thm}
\pf
For $u, v\in V$, we have
\begin{align*}
\left(\sum_{k=0}^{N}[u]_{kk}\right)\diamond 
\left(\sum_{k=0}^{N}[v]_{kk}\right)
& =\sum_{k=0}^{N}
\res_{x}T_{2k+1}((x+1)^{-k-1})(1+x)^{k}\left[Y_{V}((1+x)^{L(0)}u, x)v\right]_{kk}\nn
& =\sum_{k=0}^{N}
\res_{x}\sum_{m=0}^{k}\binom{-k-1}{m}x^{-k-m-1}
(1+x)^{k}\left[Y_{V}((1+x)^{L(0)}u, x)v\right]_{kk}\nn
& =\sum_{k=0}^{N}[u *_{k} v]_{kk}\nn
& \simeq \sum_{k=0}^{N}[u *_{N} v]_{kk} \mod Q^{\infty}(V),
\end{align*}
where in the last step, we have used the result obtained in 
the proof of Proposition 2.4 in \cite{DLM} that 
$u *_{k} v$ is equal to $u *_{N} v$ modulo $O_{k}(V)$ for $k=0, \dots, N$
and the fact that $[O_{k}(V)]_{kk}\in O^{\infty}(V)\subset Q^{\infty}(V)$.
This calculation shows that 
$A^{NN}(V)$ is closed under $\diamond$ and is thus a subalgebra
of $A^{\infty}(V)$.
Let
$f^{N}: U^{NN}(V) \to A_{N}(V)$
be defined by 
$$f^{N}\left(\sum_{k=0}^{N}[v]_{kk}\right)=v+O_{N}(V)$$ 
for $v\in V$. 

We now view $A_{N}(V)$ as an $A_{N}(V)$-module. We construct a lower-bounded generalized $V$-module 
$S(A_{N}(V))$ from $A_{N}(V)$ using  the construction in
 Section 5 of \cite{H-const-twisted-mod} as follows: Take the generating fields for 
the grading-restricted vertex algebra $V$ to be 
$Y_{V}(v, x)$ for $v\in V$. 
Take $M$ in Section 5 of \cite{H-const-twisted-mod}
to be $A_{N}(V)$.  We define the operator $L_{M}(0)$ on $M$ to be the multiplication by the scalar $N$. 
So $M$ itself is an eigenspace of $L_{M}(0)$ with eigenvalue $N$. 
Take the automorphism $g$ of $V$ in Section 5 of \cite{H-const-twisted-mod} to be $1_{V}$ 
since we are interested only in untwisted modules. Take $B$ 
in Section 5 of \cite{H-const-twisted-mod} to be $0$. Then we obtain a 
lower-bounded generalized $V$-module $\widehat{M}_{0}^{[1_{V}]}$, which shall be
denoted  by $\widehat{A_{N}(V)}$ here. 
By Theorem 3.3 in \cite{H-exist-twisted-mod} and  the construction in 
Section 5 of \cite{H-const-twisted-mod}, $\widehat{A_{N}(V)}$ is 
spanned by elements of the form 
$$(Y_{\widehat{A_{N}(V)}})_{n}(u)L_{\widehat{A_{N}(V)}}(-1)^{p}(v+O_{N}(V))$$
for homogeneous $u, v\in V$, $p\in \N$ and $n\in \wt u+N+p-1-\N$.
Let $J$ be the generalized $V$-submodule of $\widehat{A_{N}(V)}$ generated by elements of the form 
$$(Y_{\widehat{A_{N}(V)}})_{\swt u-1}(u)(v+O_{N}(V))-u*_{N}v+O_{N}(V)$$
for $u, v\in V$. Let $S(A_{N}(V))=\widehat{A_{N}(V)}/J$. Then 
$$S(A_{N}(V))=\coprod_{n\in \N} (S(A_{N}(V)))_{[n]}$$
 is a lower-bounded generalized $V$-module such that 
$(S(A_{N}(V)))_{[N]}=A_{N}(V)$. From the construction 
in Section 5 of \cite{H-const-twisted-mod}
and the definition of $S(A_{N}(V))$
above, elements of the form 
$$(Y_{S(A_{N}(V))})_{\swt u-1+N}(u)(v+O_{N}(V))$$
for homogeneous nonzero $u\in V$ and 
$v\in V\setminus O_{N}(V)$ are not $0$. Thus for 
$v\in V\setminus O_{N}(V)$, 
$v+O_{N}(V)\in A_{N}(V)$
 is not in $\Omega_{N-1}(S(A_{N}(V)))$. In other words, if $v+O_{N}(V)\in \Omega_{N-1}(S(A_{N}(V)))$,
then $v\in O_{N}(V)$. 
On the other hand, we know that $A_{N}(V)=(S(A_{N}(V)))_{[N]}\\ \subset \Omega_{N}(S(A_{N}(V)))$.

Let $W$ be a lower-bounded generalized 
$V$-module. Then $\ker \vartheta_{Gr(W)}$ is a two-sided ideal of $U^{\infty}(V)$. So 
$\ker \vartheta_{Gr(W)}\cap U_{NN}(V)$ is a 
two-sided ideal of $U_{NN}(V)$. 
From \cite{DLM}, the map $o_{W}: V\to {\rm End}\;\Omega_{N}(W)$ defined by 
$o_{W}(v)=(Y_{W})_{\swt v-1}(v)=\res_{x}x^{-1}Y_{W}(x^{L_{V}(0)}v, x)$ gives $\Omega_{N}(W)$ an
$A_{N}(V)$-module structure. In particular, $o_{W}(O_{N}(V))=0$. So $O_{N}(V)\subset \ker o_{W}$.
We take $W=S(A_{N}(V))$.  Then
$A_{N}(V)$ is an $A_{N}(V)$-submodule of
$\Omega_{N}(S(A_{N}(V)))$. We use $o_{A_{N}(V)}$ to denote the 
corresponding map from $V$ to ${\rm End}\;A_{N}(V)$. 
By the definition of $o_{A_{N}(V)}$, we have 
$$o_{A_{N}(V)}(u)(v+O_{N}(V))
=u*_{N}v+O_{N}(V)$$
for $u, v\in V$. 
For $u\in \ker o_{A_{N}(V)}$, we have
$$o_{A_{N}(V)}(u)(v+O_{N}(V))=0$$ 
for $v\in V$. So we have 
$u*_{N}v+O_{N}(V)=0$ or $u*_{N}v\in O_{N}(V)$. In particular, for $v=\one$, we have 
$u*_{N}\one\in O(V)$. But modulo $O_{N}(V)$, $u*_{N}\one$ is equal to $u$. So $u\in O_{N}(V)$. 
This means $\ker o_{S_{N}(A_{N}(V))}\subset O_{N}(V)$ and thus 
$\ker o_{A_{N}(V)}= O_{N}(V)$. 

For $v\in V$, we have shown that
$v+O_{N}(V) \in \Omega_{N-1}(S(A_{N}(V)))$ implies $v\in O_{N}(V)$
or equivalently $v+O_{N}(V)=O_{N}(V)$. So using our notation 
above, we see that 
$$[v+O_{N}(V)]_{N}=(v+O_{N}(V))+\Omega_{N-1}(S(A_{N}(V)))$$
is an element 
of $Gr_{N}(S(A_{N}(V)))$ and if it is equal to $0\in Gr_{N}(S(A_{N}(V)))$, 
then $v+O_{N}(V)=O_{N}(V)$. 
By definition, for  $u, v\in V$,
\begin{align*}
\vartheta_{Gr_{N}(S(A_{N}(V)))}([u]_{NN})[(v+O_{N}(V))]_{N}&=
[\res_{x}x^{-1}Y_{S(A_{N}(V))}(x^{L_{V}(0)}u, x)(v+O_{N}(V))]_{N}\nn
&=[o_{A_{N}(V)}(u)(v+O_{N}(V))]_{N}.
\end{align*}
Then  $\vartheta_{Gr_{N}(S(A_{N}(V)))}([u]_{NN})[(v+O_{N}(V))]_{N}
=0$ if and only if 
$o_{A_{N}(V)}(u)(v+O_{N}(V))\in \Omega_{N-1}(S(A_{N}(V)))$. 

For $\sum_{k=0}^{N}[u]_{kk}\in Q^{\infty}(V)$, 
\begin{align*}
\vartheta_{S_{N}(A_{N}(V))}([u]_{NN})
[(v+O_{N}(V))]_{N}
&=\vartheta_{S_{N}(A_{N}(V))}\left(\sum_{k=0}^{N}[u]_{kk}\right)
[(v+O_{N}(V))]_{N}\nn
&=0
\end{align*}
for all $v\in V$. So $o_{A_{N}(V)}(u)(v+O_{N}(V))\in 
\Omega_{N-1}(S(A_{N}(V)))$ or equivalently 
$[o_{A_{N}(V)}(u)(v+O_{N}(V))]_{N}$ 
is equal to $0\in 
Gr_{N}(S_{N}(A_{N}(V)))$ for all $v\in V$. Thus $o_{S_{N}(A_{N}(V))}(u)
(v+O_{N}(V))=O_{N}(V)$ or equivalently $o_{S_{N}(A_{N}(V))}(u)
(v+O_{N}(V))$ is equal to $0\in A_{N}(V)$. Then we have
$u\in \ker o_{S_{N}(A_{N}(V))}= O_{N}(V)$. 

We have proved that 
$\sum_{k=0}^{N}[u]_{kk}\in Q^{\infty}(V)$ implies 
$u\in \ker o_{S_{N}(A_{N}(V))}= O_{N}(V)$. On the other hand, 
since 
$O_{N}(V)\subset O_{k}(V)$ for $k=0, \dots, N$ and $[O_{k}(V)]_{kk}
\subset O^{\infty}(V)\subset Q^{\infty}(V)$, we have 
$\sum_{k=0}^{N}[u]_{kk}\in Q^{\infty}(V)$
for $u\in O_{N}(V)$. Thus $\sum_{k=0}^{N}[u]_{kk}\in Q^{\infty}(V)$ 
if and only if $u\in O_{N}(V)$. By this result, we obtain 
$\ker f^{N}=U^{NN}\cap Q^{\infty}(V)$. In particular, 
$f^{N}$ induces a linear isomorphism $\hat{f}^{N}: A^{NN}(V)
\to A_{N}(V)$. 

For $u, v\in V$, using the calculation above, we have
\begin{align*}
\hat{f}&^{N}\left(\left(\sum_{k=0}^{N}[u]_{kk}+Q^{\infty}(V)\right)\diamond 
\left(\sum_{k=0}^{N}[v]_{kk}+Q^{\infty}(V)\right)\right)\nn
&=\hat{f}^{N}\left(\left(\sum_{k=0}^{N}[u]_{kk}\right)\diamond 
\left(\sum_{k=0}^{N}[v]_{kk}\right)+Q^{\infty}(V)\right)\nn
&=\hat{f}^{N}\left(\sum_{k=0}^{N}[u *_{N} v]_{kk}+Q^{\infty}(V)\right)\nn
&=u*_{N} v+O_{N}(V)\nn
&=(u+O_{N}(V))*_{N} (v+O_{N}(V)).
\end{align*}
Therefore $\hat{f}^{N}$ is an  isomorphism of associative algebras. 

Since $\one+O_{N}(V)$ is the identity of $A_{N}(V)$,
$\one^{N}+O^{\infty}(V)$ is the identity of $A^{NN}(V)$. If $V$ is a vertex operator 
algebra with the conformal vector $\omega$, since
 $\omega+O_{N}(V)$ is in the center of $A_{N}(V)$, 
$\omega^{N}+O^{\infty}(V)$ is in the center of $A^{NN}(V)$. 
\epfv

\subsection{Associative algebras from finite matrices}

We now introduce new subalgebras  of $A^{\infty}(V)$. For $N\in \N$, 
let $U^{N}(V)$ be the space of all $(N+1)\times (N+1)$ 
matrices with entries in $V$. It is clear that $U^{N}(V)$ can be canonically
embedded into $U^{\infty}_{0}(V)$ as a subspace. We shall view $U^{N}(V)$ 
as a subspace of $U^{\infty}_{0}(V)$ in this paper. 
As a subspace of   $U^{\infty}_{0}(V)$, $U^{N}(V)$ 
consists of infinite matrices in $U^{\infty}(V)$
whose $(k, l)$-th entries for $k> N$ or $l> N$ are all $0$ and is
 spanned by elements of the form 
$[v]_{kl}$ for $v\in V$, $k, l=0, \dots, N$.

Recall the element 
$$\one^{N}=\sum_{k=0}^{N}[\one]_{kk},$$
that is, $\one^{N}$ is the element of $U^{N}(V)$ with the only nonzero entries to be 
equal to $\one$ at the diagonal $(k, k)$-th entries for $k=0, \dots, N$. 
By (\ref{defn-diamond}), we have 
$$\one^{N}\diamond [v]_{kl}=
\res_{x}T_{k+l+1}((x+1)^{-l-1})
(1+x)^{l}[Y_{V}((1+x)^{L(0)}\one, x)v]_{kl}=[v]_{kl}$$
for $v\in V$ and $k, l=0, \dots, N$. So 
$\one^{N}$ is a left identity of 
$U^{N}(V)$ with respect to the product $\diamond$. 
Note that for $v\in V$ and $k, l=0, \dots, N$, 
$$[v]_{kl}\diamond\one^{N}=\res_{x}T_{k+l+1}((x+1)^{-k-1})
(1+x)^{l}[Y_{V}((1+x)^{L(0)}v, x)\one]_{kl}=[v]_{kl}\diamond\one^{\infty}.$$
This formula together with (\ref{one-on-right}) immediately gives 
\begin{align}\label{one-N-on-right}
[v]_{kl}\diamond\one^{N}
&= \sum_{m=0}^{k}\binom{-k-1}{m}\left[\binom{L_{V}(-1)+L_{V}(0)+l}{k+m}v\right]_{kl}
\end{align}
for $v\in V$ and $k, l=0, \dots, N$.

By  (\ref{defn-diamond-1}), for $u, v\in V$ and $k, n, l=0, \dots, N$,
\begin{equation}\label{defn-diamond-1-N}
[u]_{kn}\diamond [v]_{nl}
=\res_{x}T_{k+l+1}((x+1)^{-k+n-l-1})(1+x)^{l}\left[Y_{V}((1+x)^{L(0)}u, x)v\right]_{kl}\in U^{N}(V).
\end{equation}
So $U^{N}(V)$ is closed under the product $\diamond$. 
Let
$$A^{N}(V)=\{\mathfrak{v}+Q^{\infty}(V)\;|\;\mathfrak{v}\in U^{N}(V)\}
=\pi_{A^{\infty}(V)}(U^{N}(V)),$$
where 
$\pi_{A^{\infty}(V)}$ is the projection from $U^{\infty}(V)$ to 
$A^{\infty}(V)$. Then $A^{N}(V)$ is
spanned by elements of the form $[v]_{kl}+Q^{\infty}(V)$ for $v\in V$ and $k, l=0, \dots, N$.

\begin{prop}
The subspace $A^{N}(V)$ is closed under $\diamond$ and is thus a  subalgebra
of $A^{\infty}(V)$ with the identity $\one^{N}+Q^{\infty}(V)$. 
\end{prop}
\pf
By (\ref{defn-diamond-1-N}),
we have 
\begin{align*}
&([u]_{kn}+Q^{\infty}(V))\diamond ([v]_{nl}+Q^{\infty}(V))\nn
&\quad =\res_{x}T_{k+l+1}((x+1)^{-k+n-l-1})
(1+x)^{l}\left[Y_{V}((1+x)^{L(0)}u, x)v\right]_{kl}+Q^{\infty}(V)\nn
&\quad \in 
A^{N}(V)
\end{align*}
for $u, v\in V$ and $k, n, l=0, \dots, N$. 
Thus $A^{N}(V)$  is closed under $\diamond$ and is thus a subalgebra
of $A^{\infty}(V)$. 

Since $\one^{N}$ is a left identity of 
$U^{N}(V)$ with respect to the product $\diamond$, $\one^{N}+Q^{\infty}(V)$
is a left identity of $A^{N}(V)$. Since
$$[v]_{kl}\diamond \one^{N}=[v]_{kl}
\diamond \one^{\infty}\equiv [v]_{kl}\mod Q^{\infty}(V),$$
$\one^{N}+Q^{\infty}(V)$ is also a right identity of $A^{N}(V)$. In particular, it is the identity of $A^{N}(V)$.
\epfv

\begin{rema}
{\rm We have derived $A^{N}(V)$ as a subalgebra of $A^{\infty}(V)$. One can certainly 
obtain $A^{N}(V)$ directly starting with the space $U^{N}(V)$ of $(N+1)\times (N+1)$ matrices 
with entries in $V$. }
\end{rema}

\begin{rema}
{\rm It is clear from the definition that $A^{nn}(V)$ for $n=0, \dots, N$ are subalgebras of $A^{N}(V)$. In
particular, the Zhu algebra $A(V)$ in \cite{Z} and its generalizations $A_{n}(V)$ for $n=0, \dots, N$ by Dong, Li and Mason 
in \cite{DLM} can be viewed as subalgebras of $A^{N}(V)$. In the case $N=0$, $A^{0}$ is equal 
to $A^{00}=A_{00}(V)$ and is thus isomorphic to 
the Zhu algebra $A(V)$ by Theorem \ref{zhu-isom}. }
\end{rema}

We say that $V$ is of positive energy
if $V=\coprod_{n\in \N}V_{(n)}$ and $V_{(0)}=\C\one$. (In some papers, $V$ being of positive energy is 
said to be of CFT type.) We recall that for $n\in \N$, $V$ is $C_{n}$-cofinite if $\dim V/C_{n}(V)<\infty$, where $C_{n}(V)$
is the subspace of $V$ spanned by elements of the form $(Y_{V})_{-n}(u)v$ for $u, v\in V$.

\begin{thm}\label{finite-d}
Assume that  $V$ is of positive energy and $C_{2}$-cofinite. Then $A^{N}(V)$ is finite dimensional. 
\end{thm}
\pf
By Theorem 11 in 
\cite{GN} (see Proposition 5.5 in \cite{AN}), 
$V$ is also $C_{n}$-cofinite for $n\ge 2$. In particular, $V$ is $C_{k+l+2}$-cofinite
for $k, l=0, \dots, N$. 
By definition, $C_{k+l+2}(V)$ are spanned by elements of the 
form $(Y_{V})_{-k-l-2}(u)v$ for $u, v\in V$. Since $V$ is $C_{k+l+2}$-cofinite,
there exists a finite dimensional subspace $X_{k+ l}$ of $V$ such that 
$X+C_{k+l+2}(V)=V$. Let $U^{N}(X)$ be the subspace of $U^{N}(V)$ consisting 
matrices in $U^{N}(V)$ whose $(k, l)$-th 
entries are in $X$ for $k, l=0, \dots, N$. Since $X_{k+l}$ for 
$k, l=0, \dots, N$ are finite dimensional, 
$U^{N}(X)$ is also finite dimensional. We now prove $U^{N}(X)+(O^{\infty}(V)\cap U^{N}(V))=U^{N}(V)$. 
To prove this, we need only prove that every element of $U^{N}(V)$ of the form 
$[v]_{kl}$ for $v\in V$ and $0\le k, l\le N$, can be written as $[v]_{kl}=[v_{1}]_{kl}+[v_{2}]_{kl}$, where $v_{1}\in X_{k+l}$ and 
$v_{2}\in V$ such that $[v_{2}]_{kl}\in O^{\infty}(V)$. We shall denote the subspace of $V$ consisting of elements
$v$ such that $[v]_{kl}\in O^{\infty}(V)$ by $O^{\infty}_{kl}(V)$. 
Then what we need to prove is $V=X_{k+l}+O^{\infty}_{kl}(V)$.

We can always take $X_{k+l}$ to be a subspace of $V$ containing $\one$. 
We use
induction on the weight of $v$. When $\wt v=0$, $v$ is proportional to $\one$  and
can indeed be written as $v=v+0$, where $v\in X$ and $0\in O^{\infty}_{kl}(V)$.

Assume that  when $\wt v=p<q$,  $v=v_{1}+v_{2}$, where  $v_{1}\in X_{k+l}$ and 
$v_{2}\in O_{kl}^{\infty}(V)$. 
Then since $V$ is $C_{k+l+2}$-cofinite, for $v\in V_{(q)}$, there exists
homogeneous $u_{1}\in X_{k+l}$ and homogeneous $u^{i}, v^{i}\in V$ for $i=1, 
\dots, m$ such that $v=u_{1}+\sum_{i=1}^{m}u^{i}_{-k-l-2}v^{i}$.
Moreover, we can always find such $u_{1}$ and $u^{i}, v^{i}\in V$ for $i=1, 
\dots, m$ such that $\wt u_{1}=\wt u^{i}_{-k-l-2}v^{i}=\wt v=q$. 
Since 
$$\wt u^{i}_{n-k-l-2}v^{i}<\wt u^{i}_{-k-l-2}v^{i}=\wt v=q$$
for $i=1, 
\dots, m$ and
$n\in \Z_{+}$, by induction assumption, 
$u^{i}_{n-k-l-2}v^{i}\in X_{k+l}+O^{\infty}_{kl}(V)$
for $i=1, \dots, m$ and $k\in \Z_{+}$.
Thus 
\begin{align*}
v&=u_{1}+\sum_{i=1}^{m}u^{i}_{-k-l-2}v^{i}\nn
&=u_{1}+\sum_{i=1}^{m}
\res_{x}x^{-k-l-2}(1+x)^{l}Y((1+x)^{L(0)}u^{i}, x)v^{i}\nn
& \quad-\sum_{i=1}^{m}\sum_{n\in \Z_{+}}\binom{\wt u^{i}+l}{n}
u^{i}_{n-k-l-2}v^{i}.
\end{align*}
By definition, 
$$[\res_{x}x^{-k-l-2}(1+x)^{l}Y((1+x)^{L(0)}u^{i}, x)v^{i}]_{kl}\in O^{\infty}(V).$$
Thus 
$$\res_{x}x^{-k-l-2}(1+x)^{l}Y((1+x)^{L(0)}u^{i}, x)v^{i}\in O^{\infty}_{kl}(V).$$
Thus we have $v=v_{1}+v_{2}$, where  $v_{1}\in X_{k+l}$ and 
$v_{2}\in O_{kl}^{\infty}(V)$. 
By induction principle, we have $V=X_{k+l}+O^{\infty}_{kl}(V)$. 

We now have proved $U^{N}(X)+(O^{\infty}(V)\cap U^{N}(V))=U^{N}(V)$. 
Since $O^{\infty}(V)\cap U^{N}(V)\subset Q^{\infty}(V)\cap U^{N}(V)$, we also have 
$U^{N}(X)+(Q^{\infty}(V)\cap U^{N}(V))=U^{N}(V)$. Since $U^{N}(X)$ is finite
dimensional, $A^{N}(V)$ is finite dimensional. 
\epfv

\setcounter{equation}{0} \setcounter{thm}{0} 

\section{Lower-bounded generalized $V$-modules and graded $A^{N}(V)$-modules}

By Theorem \ref{assoc-alg},  the associated graded space $Gr(W)$ of a filtration of a lower-bounded 
generalized $V$-module $W$  is a nondegenerate graded $A^{\infty}(V)$-module.  In this section, for $N\in \N$, we 
give an $A^{N}(V)$-module structure to a subspace of $Gr(W)$ and use it
to study $W$. 

Let $N\in \N$. Let $W$ be a lower-bounded generalized $V$-module.
Since $A^{N}(V)$ is a subalgebra of $A^{\infty}(V)$, $Gr(W)$ as an $A^{\infty}(V)$-module 
is also an $A^{N}(V)$-module. 
Let 
$$Gr^{N}(W)=\coprod_{n=0}^{N}Gr_{n}(W)\subset Gr(W).$$
By the definition of $\vartheta_{Gr(W)}$, we see that for $\mathfrak{v}\in A^{N}(V)$ and $[w]_{n}\in Gr^{N}(W)$,
$\vartheta_{Gr(W)}(\mathfrak{v})[w]_{n}\in Gr^{N}(W)$. Thus 
$Gr^{N}(W)$  is an $A^{N}(W)$-submodule of  $Gr(W)$.
But $Gr^{N}(W)$ has some additional structures and properties and we are only interested
in those $A^{N}(W)$-modules having these additional structures and properties. 
Similar to Definition \ref{N-graded-A-inf-mod}, we have the following notion:

\begin{defn}\label{gr-A-N-mod}
{\rm Let $M$ be an $A^{N}(V)$-module $M$ with the $A^{N}(V)$-module structure on $M$ given by 
$\vartheta_{M}: A^{N}(V)\to {\rm End}\; M$. We say that $M$ is a {\it graded $A^{N}(V)$-module}
if the following conditions are 
satisfied:
\begin{enumerate}

\item  $M=\coprod_{n=0}^{N}G_{n}(M)$ such that for $v\in V$ and $k, l=0, \dots, N$,
$\vartheta_{M}([v]_{kl}+Q^{\infty}(V))$ maps $G_{n}(M)$ for $0\le n\le N$ to 
$0$ when $n\ne l$ and to $G_{k}(M)$ when $n=l$.

\item $M$ is  a direct sum of generalized eigenspaces of of an operator $L_{M}(0)$
on $M$. $G_{n}(M)$ for $n\in \N$ are invariant under $L_{M}(0)$ and 
the real parts of the eigenvalues of $L_{M}(0)$ has a lower bound. 

\item There is a linear map $L_{M}(-1): \coprod_{n=0}^{N-1}G_{n-1}(M)\to 
\coprod_{n=1}^{N}G_{n}(M)$
mapping $G_{n}(M)$ to $G_{n+1}(M)$ for $n=0, \dots, N-1$.

\item The commutator relations
\begin{align*}
{[L_{M}(0), L_{M}(-1)]}&=L_{M}(-1),\nn
{[L_{M}(0), \vartheta_{M}([v]_{kl}+Q^{\infty}(V))]}&=(k-l)\vartheta_{M}([v]_{kl}+Q^{\infty}(V)),\nn
{[L_{M}(-1), \vartheta_{M}([v]_{pl}+Q^{\infty}(V))]}&=\vartheta_{M}([L_{V}(-1)v]_{(p+1)l}+Q^{\infty}(V))
\end{align*}
hold for $v\in V$, 
$k, l=0, \dots, N$ and $p=0, \dots, N-1$. 

\end{enumerate}
A graded $A^{N}(V)$-module $M$ is said to be {\it nondegenerate} if 
the following additional condition holds: For $w\in G_{l}(M)$, if $\vartheta_{M}([v]_{0l}+Q^{\infty}(V))w=0$ 
for all $v\in V$, then $w=0$. 
Let $M_{1}$ and $M_{2}$ be graded $A^{N}(V)$-modules.
An {\it graded $A^{N}(V)$-module map} from 
$M_{2}$ to $M_{2}$ is an $A^{N}(V)$-module map $f: M_{1}\to M_{2}$ 
such that $f(G_{n}(M_{1}))\subset G_{n}(M_{2})$ for $n=0 \dots, N$, 
$f\circ L_{M_{1}}(0)=L_{M_{2}}(0)\circ f$ and $f\circ L_{M_{1}}(-1)=L_{M_{2}}(-1)\circ f$.
A {\it graded $A^{N}(V)$-submodule} of a graded $A^{N}(V)$-module $M$ is an $A^{N}(V)$-submodule 
$M_{0}$ of $M$ such that  with the $A^{N}(V)$-module structure, the $\N$-grading induced from $M$
and the operators $L_{M}(0)\mbar_{M_{0}}$ and $L_{M}(-1)\mbar_{M_{0}}$, 
$M_{0}$ is a graded $A^{N}(V)$-module.
A graded $A^{\infty}(V)$-module $M$ is said to be {\it generated by a subset $S$}  if 
$M$  is equal  to the smallest graded $A^{N}(V)$-submodule containing $S$, or equivalently, 
$M$  is spanned by homogeneous elements 
obtained by applying elements of $A^{N}(V)$, $L_{M}(0)$ and $L_{M}(-1)$ to 
homogeneous summands of  elements of $S$.  A graded $A^{N}(V)$-module is said to be 
{\it irreducible} if it has no nonzero proper graded $A^{N}(V)$-modules. A graded $A^{N}(V)$-module 
is said to be {\it completely reducible} if it is a direct sum of irreducible graded $A^{N}(V)$-modules.}
\end{defn}

From the discussion above and 
the property of $Gr^{N}(W)$, we obtain immediately:

\begin{prop}\label{A-N-V-mod}
For a lower-bounded generalized $V$-module $W$, $Gr^{N}(W)$
is a nondegenerate graded $A^{N}(V)$-module.  Let $W_{1}$ and $W_{2}$ be lower-bounded generalized $V$-modules 
and $f: W_{1}\to W_{2}$ a $V$-module map. Then $f$ induces a graded $A^{N}(V)$-module map
$Gr^{N}(f): Gr^{N}(W_{1})\to Gr^{N}(W_{2})$.
\end{prop}

We have the following results on irreducible and completely reducible
lower-bounded generalized $V$-modules without additional conditions:

\begin{prop}\label{T-N-preserv-irr}
Let $W$ be a lower-bounded generalized $V$-module. 
If $W$ is irreducible or completely reducible, then $Gr^{N}(W)$ is equivalent to $T_{N}(W)$ as an $A^{N}(V)$-module and 
is also irreducible or completely reducible, respectively.
\end{prop}
\pf 
Let $W$ be irreducible. By Proposition \ref{W=GrW}, 
$\Omega_{n}(W)=T_{n}(W)$ for $n=0, \dots, N$. Then $T_{N}(W)$ is a nondegenerate graded $A^{N}(V)$-module 
equivalent to $Gr^{N}(W)$. 
We need to prove that the nondegenerate  graded $A^{N}(V)$-module $T_{N}(W)$ is 
irreducible. 

Let $M$ be a nonzero graded
$A^{N}(V)$-submodule of $T_{N}(W)$. 
We use  the construction in Section 5 of \cite{H-const-twisted-mod} to construct 
a universal lower-bounded generalized $V$-module $\widehat{M}$ from $M$.  
We take the generating fields for the grading-restricted vertex algebra $V$ to be 
$Y_{V}(v, x)$ for $v\in V$. By definition, 
$M$ is a direct sum of generalized eigenspaces of $L_{M}(0)$ and the real parts of the eigenvalues of 
$L_{M}(0)$ have a lower bound $B\in \R$. We take $M$ and $B$ in Section 5 of \cite{H-const-twisted-mod}
to be the given nondegenerate  graded $A^{N}(V)$-module $M$ and the lower bound $B$ above. 
Using the construction in Section 5 of \cite{H-const-twisted-mod},
we obtain a universal lower-bounded generalized $V$-module 
$\widehat{M}^{[1_{V}]}_{B}$. For simplicity, we shall denote it simply by 
$\widehat{M}$. By the universal property of $\widehat{M}$ (Theorem 5.2 in \cite{H-const-twisted-mod}), 
for the embedding map 
$e_{M}: M\to T_{N}(W)$, 
there is a unique $V$-module map 
$\widehat{e_{M}}: \widehat{M}\to W$ such that $\widehat{e_{M}}\mbar_{M}=e$.
Then $\widehat{e_{M}}(\widehat{M})$  is a generalized $V$-submodule of $W$
generated by $M$. 
It is nonzero since $M\subset 
\widehat{e_{M}}(\widehat{M})$. Since $W$ is irreducible, it must be $W$. Then $W$ is generated by $M$. 
In particular, $T_{N}(W)$ is obtained by applying 
the components of the vertex operators on $W$, $L_{W}(0)$ and $L_{W}(-1)$  to elements of $M$.  Since the components of 
the vertex operators on $W$  and the operators $L_{W}(0)$ and 
 $L_{W}(-1)$ preserving $T_{N}(W)$
are by definition the actions of elements of $A^{N}(V)$, $L_{W}(0)$ and $L_{M}(-1)$ preserving $T_{N}(W)$, 
we see that as a graded $A^{N}(V)$-module, 
$T_{N}(W)$ is generated by $M$. But $M$ itself is an $A^{N}(V)$-submodule of $T_{N}(W)$.
So we have $M=T_{N}(W)$. Thus $T_{N}(W)$ as a graded  $A^{N}$-module is 
irreducible. 

If $W$ is completely reducible, by Proposition \ref{W=GrW} again, $\Omega_{n}(W)=T_{n}(W)$ for $n=0, \dots, N$.
Then $T_{N}(W)$ is a nondegenerate  graded $A^{N}(V)$-module 
equivalent to $Gr^{N}(W)$. Since $W$ is completely reducible, $W=\coprod_{\mu\in \mathcal{M}}W^{\mu}$,
where $W^{\mu}$ for $\mu\in \mathcal{M}$ are irreducible lower-bounded generalized $V$-modules. 
By the definition of $T_{N}(W)$, we have $T_{N}(W)=\coprod_{\mu\in\mathcal{M}}T_{N}(W^{\mu})$.
From what we have proved above, for $\mu\in \mathcal{M}$, 
$T_{N}(W^{\mu})$ is an irreducible graded $A^{N}(V)$-module. Thus we see that $T_{N}(W)$ is completely reducible.
\epfv

Let $M$ be a 
graded $A^{N}(V)$-module given by a linear map $\vartheta_{M}: A^{N}(V)\to {\rm End}\;M$ and operators 
$L_{M}(0)$ and $L_{M}(-1)$.
We now construct a lower-bounded generalized $V$-module $S^{N}(M)$ from $M$. 
We use  the construction in Section 5 of \cite{H-const-twisted-mod}. 
We take the generating fields for the grading-restricted vertex algebra $V$ to be 
$Y_{V}(v, x)$ for $v\in V$. By definition, 
$M$ is a direct sum of generalized eigenspaces of $L_{M}(0)$ and the real parts of the eigenvalues of 
$L_{M}(0)$ has a lower bound $B\in \R$. We take $M$ and $B$ in Section 5 of \cite{H-const-twisted-mod}
to be the given graded $A^{N}(V)$-module $M$ and the lower bound $B$ above. 
Using the construction in Section 5 of \cite{H-const-twisted-mod},
we obtain a universal lower-bounded generalized $V$-module 
$\widehat{M}^{[1_{V}]}_{B}$. For simplicity, we shall denote it simply by 
$\widehat{M}$. 

By Theorem 3.4 in \cite{H-exist-twisted-mod} and the construction in  Section 5 of \cite{H-const-twisted-mod}
and by identifying elements of the form 
$(\psi_{\widehat{M}}^{a})_{-1, 0}$ with basis elements $w^{a}\in M$ for $a\in A$
for a basis $\{w^{a}\}_{a\in A}$ of $M$, 
we see that $\widehat{M}$ is generated by $M$  (in the sense of Definition 3.1 in \cite{H-exist-twisted-mod}). 
Moreover, Theorems 3.3 and 3.4 in \cite{H-exist-twisted-mod} 
state that elements of the form 
$L_{\widehat{M}}(-1)^{p}w^{a}$ for $p\in \N$ and $a\in A$ are linearly independent and 
$\widehat{M}$ is spanned by elements obtained by applying the components of the vertex operators
to these elements. In particular, we identify
$M$ as a subspace of $\widehat{M}$. Let $J_{M}$  be  the generalized $V$-submodule 
of $\widehat{M}$ generated by elements of the forms
\begin{equation}\label{J-M-0}
\res_{x}x^{l-k-1}Y_{\widehat{M}}(x^{L_{V}(0)}v, x)w
\end{equation}
for $ l=0, \dots, N$, $k\in -\Z_{+}$ and $w\in G_{l}(M)$,
\begin{equation}\label{J-M}
\res_{x}x^{l-k-1}Y_{\widehat{M}}(x^{L_{V}(0)}v, x)w-\vartheta_{M}([v]_{kl})w
\end{equation}
for $v\in V$, $k, l=0, \dots, N$ and $w\in G_{l}(M)$ and
\begin{equation}\label{J-M-2}
L_{\widehat{M}}(-1)w-L_{M}(-1)w
\end{equation}
for $w\in \coprod_{n=0}^{N-1}G_{n}(M)$. 

Let $S^{N}(M)=\widehat{M}/J_{M}$. Then $S^{N}(M)$ is a lower-bounded generalized $V$-module. 
Let $\pi_{S^{N}(M)}$ be the projection from $\widehat{M}$ to $S^{N}(M)$. 
Since $\widehat{M}$ is generated by $M$ (in the sense of Definition 3.1 in \cite{H-exist-twisted-mod}),
$S^{N}(M)$ is generated $\pi_{S^{N}(M)}(M)$ (in the same sense). In particular, 
$S^{N}(M)$ is spanned by elements of the form 
\begin{equation}\label{S-N-M}
\res_{x}x^{(l+p)-n-1}Y_{S^{N}(M)}(x^{L_{V}(0)}v, x)L_{S^{N}(M)}(-1)^{p}\pi_{S^{N}(M)}(w)
\end{equation}
for $v\in V$, $l=0, \dots, N$, $n, p\in \N$ and $w\in G_{l}(M)$. For $n\in \N$, let $G_{n}(S^{N}(M))$ be the subspace of 
$S^{N}(M)$ spanned by elements of the form (\ref{S-N-M}) for $v\in V$, 
 $l=0, \dots, N$, $p\in \N$ and $w\in G_{l}(M)$.

\begin{prop}\label{S-M-prop}
Let $M$ be a  graded $A^{N}(V)$-module. 
\begin{enumerate}

\item For $0\le n\le  N$, $G_{n}(S^{N}(M))=\pi_{S^{N}(M)}(G_{n}(M))$ and 
for $0\le n_{1},  n_{2}\le N$, $n_{1}\ne n_{2}$, $G_{n_{1}}(S^{N}(M))\cap G_{n_{2}}(S^{N}(M))=0$. 
Moreover, $S^{N}(M)=\coprod_{n\in \N}G_{n}(S^{N}(M))$ and $\pi_{S^{N}(M)}(M)
=\coprod_{n=0}^{N}G_{n}(S^{N}(M))$. 

\item For $n=0, \dots, N$, 
\begin{equation}\label{S-M-prop-1}
\pi_{S^{N}(M)}\left(\coprod_{j=0}^{n}G_{j}(M)\right)
=\coprod_{j=0}^{n}G_{j}(S^{N}(M))\subset  
\Omega_{n}(S^{N}(M))
\end{equation}
and in the case that $M$ is nondegenerate,
\begin{equation}\label{S-M-prop-2}
\pi_{S^{N}(M)}\left(\coprod_{j=n}^{N}G_{j}(M)\right)\cap \Omega_{n}(S^{N}(M))
=\left(\coprod_{j=n}^{N}G_{j}(S^{N}(M))\right)\cap \Omega_{n}(S^{N}(M))=0. 
\end{equation}
 
\item  In the case that $M$ is nondegenerate,
$M$ is equivalent to a graded $A^{N}(V)$-submodule of $Gr^{N}(S^{N}(M))$.

\end{enumerate}
\end{prop}
\pf
By definition, $G_{n}(S^{N}(M))$ for  $0\le n\le N$ is spanned by elements of the form 
(\ref{S-N-M}) for $v\in V$, $l=0, \dots, N$, $p\in \N$ and $w\in G_{l}(M)$. 
Using the $L(-1)$-commutator formula for the vertex operator map $Y_{S^{N}(M)}$, 
we see that it is also spanned by elements of the form 
\begin{equation}\label{S-N-M-1}
\res_{x}x^{l-k-1}L_{S^{N}(M)}(-1)^{p}Y_{S^{N}(M)}(x^{L_{V}(0)}v, x)\pi_{S^{N}(M)}(w)
\end{equation}
for $v\in V$, $l, k=0, \dots, N$, $p=0, \dots, n-k$ and $w\in G_{l}(M)$. 
Since elements of the forms (\ref{J-M}) and (\ref{J-M-2}) are in $J_{M}$, we see that
(\ref{S-N-M-1}) is in fact equal to 
\begin{equation}\label{S-N-M-2}
\pi_{S^{N}(M)}(L_{M}(-1)^{p}\vartheta_{M}([v]_{kl}+Q^{\infty}(V))w)\in \pi_{S^{N}(M)}(G_{n}(M)).
\end{equation}
Since $L_{M}(-1)^{p}\vartheta_{M}([v]_{k(l+p)}+Q^{\infty}(V))w$ 
for $v\in V$,   $l, k=0, \dots, N$, $p=0, \dots, n-k$ and $w\in G_{l}(M)$ certainly   
span $G_{n}(M)$ (in fact, we need only $v=\one$, $k=l=n$, $p=0$ and $w\in G_{n}(W)$)
and elements of the form (\ref{S-N-M-1}) for $v\in V$, 
$l, k=0, \dots, N$, $p=0, \dots, n-k$ and $w\in G_{l}(M)$ span  $G_{n}(S(G))$ for 
$0\le n\le N$, we see that 
elements of the form (\ref{S-N-M-2})  for $v\in V$, $l, k=0, \dots, N$, $p=0, \dots, n-k$ and $w\in G_{l}(M)$
also  span $G_{n}(S^{N}(M))$. Thus we obtain $G_{n}(S^{N}(M))=\pi_{S^{N}(M)}(G_{n}(M))$
for $n=0, \dots, N$. 
When $n_{1}\ne n_{2}$, we know $G_{n_{1}}(M)\cap G_{n_{2}}(M)=0$.
Then $G_{n_{1}}(S^{N}(M))\cap G_{n_{2}}(S^{N}(M))=\pi_{S(G)}(G_{n_{1}}(M)\cap G_{n_{2}}(M))=0$.
Since  $S^{N}(M)$ is spanned by elements of the form (\ref{S-N-M})
for $v\in V$, $l=0, \dots, N$, $n, p\in \N$ and $w\in G_{l}(M)$, 
by the definition of $G_{n}(S^{N}(M))$, we have  $S^{N}(M)=\coprod_{n\in \N}G_{n}(S^{N}(M))$. 
Since $G_{n}(S^{N}(M))=\pi_{S^{N}(M)}(G_{n}(M))$ for $n=0, \dots, N$, we have
$$\pi_{S^{N}(M)}(M)=\coprod_{n=0}^{N}\pi_{S^{N}(M)}(G_{n}(M))=\coprod_{n=0}^{N}G_{n}(S^{N}(M)).$$

By definition, for $0\le j\le n\le N$, $G_{j}(S^{N}(M))\subset \Omega_{n}(S^{N}(M))$. 
Then for $j=0, \dots, n$, 
$$\pi_{S^{N}(M)}(G_{j}(M))=G_{j}(S^{N}(M))\subset  \Omega_{j}(S^{N}(M))\subset \Omega_{n}(S^{N}(M)).$$
So we obtain (\ref{S-M-prop-1}).
By the nondegeneracy of $M$, nonzero  elements of 
$G_{j}(M)$ for $N\ge j>n$ are not in $\Omega_{n}(\widehat{M})$. 
From the construction of $\widehat{M}$, nonzero elements of the form (\ref{J-M-0}),  (\ref{J-M})
or  (\ref{J-M-2})  are not in $M\subset \widehat{M}$.
In particular, we see that the intersection of $J(M)$ with $M$ is $0$.
So $\pi_{S^{N}(M)}\mbar_{M}$ is injective.
Since $\pi_{S^{N}(M)}\mbar_{M}$ is injective, 
we see that nonzero  elements of 
$G_{j}(S^{N}(M))=\pi_{S^{N}(M)}(G_{j}(M))$ for $N\ge j>n$ are not in $\Omega_{n}(S^{N}(M))$. 
Thus we obtain (\ref{S-M-prop-2}).

For $0\le n\le N$ and $w\in G_{n}(M)$, we define $f_{M}(w)=\pi_{S^{N}(M)}(w)+\Omega_{n-1}(S^{N}(M))$. 
Since $\pi_{S^{N}(M)}(w)\in \Omega_{n}(S^{N}(M))$, $f_{M}(w)\in Gr_{n}(S^{N}(M))$. Therefore 
we obtain a linear map $f_{M}: M\to Gr^{N}(S^{N}(M))$. 
It is clear from the definition that $f_{M}$ is in fact a  graded $A^{N}(V)$-module map. 
If for some $0\le n\le N$ and $w\in G_{n}(M)$, $f_{M}(w)=0$, then $\pi_{S^{N}(M)}(w)\in \Omega_{n-1}(S^{N}(M))$. 
But we have proved above that nonzero elements of 
$\pi_{S^{N}(M)}(G_{n}(M))$  are not in $\Omega_{n-1}(S^{N}(M))$. So $\pi_{S^{N}(M)}(w)=0$.
Since $\pi_{S^{N}(M)}\mbar_{M}$ is injective, we obtain $w=0$. So $f_{M}$ is injective. 
Thus $M$ is equivalent to the nondegenerate  graded $A^{N}(V)$-submodule $f_{M}(M)$ of $Gr^{N}(S^{N}(M))$. 
\epfv

\begin{rema}
{\rm As in the case of $S(G)$ in Section 3, our construction of the lower-bounded generalized $V$-module 
$\widehat{M}$ depends on the lower bound $B$ of the real parts of the 
eigenvalues of $L_{M}(0)$. But by Proposition \ref{S-M-prop}, $S^{N}(M)$ depends only on 
$M$, not on $B$.}
\end{rema}

\begin{thm}\label{irred-m-bij}
For $N\in \N$, the set of the equivalence classes of irreducible lower-bounded generalized $V$-modules is in bijection with 
the set of  the equivalence classes of irreducible nondegenerate graded $A^{N}(V)$-modules. 
\end{thm}
\pf
Recall the set $[\mathfrak{W}]_{\rm irr}$
of the equivalence classes of irreducible lower-bounded generalized 
$V$-modules in the proof of Theorem \ref{irred-A-inf-m-bij}. 
Let $[\mathfrak{M}^{N}]_{\rm irr}$ be  the 
set of  the equivalence classes of irreducible nondegenerate graded $A^{N}(V)$-modules.
Given an irreducible lower-bounded generalized $V$-module $W$, by Theorem \ref{T-N-preserv-irr}, 
$Gr^{N}(W)=T_{N}(W)$ is an irreducible nondegenerate graded $A^{N}(V)$-module. Thus we obtain a map $f:
[\mathfrak{W}]_{\rm irr}\to [\mathfrak{M}^{N}]_{\rm irr}$ 
given by $f([W])=[T_{N}(W)]$, where $[W]\in [\mathfrak{W}]_{\rm irr}$ is the equivalence class containing
the irreducible  lower-bounded generalized $V$-module $W$ and $[T_{N}(W)]\in [\mathfrak{M}^{N}]_{\rm irr}$ 
is the equivalence 
class containing the irreducible nondegenerate  graded $A^{N}(V)$-module $T_{N}(W)$. 

Given an irreducible nondegenerate graded $A^{N}(V)$-module $M$,
we have the lower-bounded generalized $V$-module $S^{N}(M)$ generated by $\pi_{S^{N}(M)}(M)$.
The main difference of the proof here and the the proof of Theorem \ref{irred-A-inf-m-bij} 
is that we do not know whether $S^{N}(M)$  is  irreducible. 
So we need to take a quotient of $S^{N}(M)$. 
Since $M$ is an irreducible nondegenerate graded $A^{N}(V)$-module, 
it is generated by any nonzero element.
Since $S^{N}(M)$ is generated by $\pi_{S^{N}(M)}(M)$, it is also generated by  any element
$w_{0}\in \pi_{S^{N}(M)}(M)$. 
Then by Theorem 4.7 in 
\cite{H-exist-twisted-mod}, there is a maximal generalized $V$-submodule $J_{\pi_{S^{N}(M)}(M), w_{0}}$ 
of $S^{N}(M)$ such 
that $J_{\pi_{S^{N}(M)}(M), w_{0}}$ does not contain $w_{0}$ and 
$S^{N}(M)/J_{\pi_{S^{N}(M)}(M), w_{0}}$ is irreducible.
The maximal generalized $V$-submodule $J_{\pi_{S^{N}(M)}(M), w_{0}}$ is in fact independent of 
$w_{0}\in \pi_{S^{N}(M)}(M)$. We prove this fact by proving that no nonzero element of 
$\pi_{S^{N}(M)}(M)$ is in $J_{\pi_{S^{N}(M)}(M), w_{0}}$. 
In fact, if a nonzero $w\in \pi_{S^{N}(M)}(M)$ is also in $J_{\pi_{S^{N}(M)}(M), w_{0}}$,
since the actions of components of vertex operators on $w$ are equal to the actions of elements of 
$A^{N}(V)$ and $M$ is generated also by $w$, we see that $w_{0}$ must also be in $J_{\pi_{S^{N}(M)}(M), w_{0}}$.
Contradiction. Thus $J_{\pi_{S^{N}(M)}(M), w_{0}}$ is in fact the maximal 
generalized $V$-submodule of $S^{N}(M)$ such 
that it does not contain nonzero elements of $M$. We denote it by $\widetilde{J}_{M}$, which depends only on 
$\pi_{S^{N}(M)}(M)$, or equivalently, $M$. 
Thus we obtain a map 
$g: [\mathfrak{M}^{N}]_{\rm irr}\to [\mathfrak{W}]_{\rm irr}$ given by $g([M])=[S^{N}(M)/\widetilde{J}_{M}]$. 

We still need to show that the two maps above are inverses of each other. 
Let $M$ be an irreducible nondegenerate graded $A^{N}(V)$-module. Since $S^{N}(M)/\widetilde{J}_{M}$ is 
irreducible, by Proposition \ref{W=GrW}, $Gr_{N}(S^{N}(M)/\widetilde{J}_{M})$ is an 
irreducible nondegenerate graded $A^{N}(V)$-module. By Proposition \ref{S-M-prop}, $M$ is equivalent to 
a nondegenerate 
graded $A^{N}(V)$-submodule of $Gr_{N}(S^{N}(M))$. As in the proof of Proposition \ref{S-M-prop},
we denote this equivalence by $f_{M}$.
Let $\pi_{\widetilde{J}_{M}}: S^{N}(M)\to S^{N}(M)/\widetilde{J}_{M}$ be the projection map. 
Since $\widetilde{J}_{M}\cap \pi_{S^{N}(M)}(M)=0$, $\pi_{\widetilde{J}_{M}}\mbar_{\pi_{S^{N}(M)}(M)}$
is injective and in particular, is not $0$. 
The $V$-module map $\pi_{\widetilde{J}_{M}}$ induces a graded $A^{N}(V)$-module map
$Gr^{N}(\pi_{\widetilde{J}_{M}}): Gr_{N}(S^{N}(M)) \to Gr_{N}(S^{N}(M)/\widetilde{J}_{M})$.
Since  $\pi_{\widetilde{J}_{M}}\mbar_{\pi_{S^{N}(M)}(M)}$ is not $0$, 
the restriction $Gr^{N}(\pi_{\widetilde{J}_{M}})\mbar_{f_{M}(M)}$
of $Gr^{N}(\pi_{\widetilde{J}_{M}})$ to the image of $M$ under $f_{M}$ is also not $0$. 
Consider the $A^{N}(V)$-module map $Gr^{N}(\pi_{\widetilde{J}_{M}}) \circ f_{M}:
M\to Gr_{N}(S^{N}(M)/\widetilde{J}_{M})$. Since $f_{M}$ is injective and 
$Gr^{N}(\pi_{\widetilde{J}_{M}})\mbar_{f_{M}(M)}\ne 0$,
 $Gr^{N}(\pi_{\widetilde{J}_{M}}) \circ f_{M}$ is not $0$. But 
both $M$ and $Gr_{N}(S^{N}(M)/\widetilde{J}_{M})$ are irreducible. So 
$Gr^{N}(\pi_{\widetilde{J}_{M}}) \circ f_{M}$ must be an equivalence of 
graded $A^{N}(V)$-modules. Moreover, by Proposition \ref{T-N-preserv-irr},
$Gr_{N}(S^{N}(M)/\widetilde{J}_{M})$ is equivalent to $T_{N}(S^{N}(M)/\widetilde{J}_{M})$. 
So $M$ is equivalent to $T_{N}(S^{N}(M)/\widetilde{J}_{M})$.  Thus $[M]=[T_{N}(\widehat{M}/\widetilde{J}_{M})]$.
This means $f(g([M]))=[M]$. So we obtain $f\circ g=1_{[\mathfrak{M}^{N}]_{\rm irr}}$. 

Let $W$ be an irreducible  lower-bounded generalized $V$-module. By Theorem \ref{T-N-preserv-irr},
$T_{N}(W)$ is an irreducible $A^{N}(V)$-module. We then have a lower-bounded generalized  $V$-module 
$S^{N}(T_{N}(W))$. 
By the universal property of $\widehat{T_{N}(W)}$,
there is a unique $V$-module map $\widehat{1_{T_{N}(W)}}: \widehat{T_{N}(W)}\to W$ such that 
$\widehat{1_{T_{N}(W)})}\mbar_{T_{N}(W)}=1_{T_{N}(W)}$, where $1_{T_{N}(W)}$
is the identity operator on $T_{N}(W)$. 
Since $W$ is irreducible, the image of $\widehat{T_{N}(W)}$ under $\widehat{1_{T_{N}(W)}}$ 
is either $0$ or $W$. Since $\widehat{1_{T_{N}(W)})}\mbar_{T_{N}(W)}=1_{T_{N}(W)}$,
 the image of $\widehat{T_{N}(W)}$ under $\widehat{1_{T_{N}(W)}}$  cannot be $0$ and thus  must be $W$.
In particular, $\widehat{1_{T_{N}(W)}}$ is surjective. Moreover, since
$J_{T_{N}(W)}$ is generated by (\ref{J-M-0}), (\ref{J-M}) and (\ref{J-M-2}) with $M=T_{N}(W)$, 
the image of $J_{T_{N}(W)}$
under $\widehat{1_{T_{N}(W)}}$  is $0$, that is, $J_{T_{N}(W)}\in \ker  \widehat{1_{T_{N}(W)}}$. 
In particular, $\widehat{1_{T_{N}(W)}}$
induces a surjective $V$-module map $f_{T_{N}(W)}: S^{N}(T_{N}(W))=\widehat{T_{N}(W)}/J_{T_{N}(W)}\to W$. 
Since $J_{T_{N}(W)}\cap T_{N}(W)=0$, $f_{T_{N}(W)}(T_{N}(W))=T_{N}(W)$. 
We have a maximal generalized $V$-submodule 
$\widetilde{J}_{T_{N}(W)}$ of $S^{N}(T_{N}(W))$ as in the construction above
such that $T_{N}(W)\cap \widetilde{J}_{T_{N}(W)}=0$  and 
$S^{N}(T_{N}(W))/\widetilde{J}_{T_{N}(W)}$ is irreducible. Since $f_{T_{N}(W)}(T_{N}(W))=T_{N}(W)$,
$\ker f_{T_{N}(W)}$ is a generalized $V$-submodule of $S^{N}(T_{N}(W))$ that does not 
contain nonzero elements of $M$. Hence $\ker f_{T_{N}(W)}\subset \widetilde{J}_{T_{N}(W)}$.
Thus we obtain a surjective $V$-module map from $S^{N}(T_{N}(W))/\widetilde{J}_{T_{N}(W)}$ to $W$.
Since both $S^{N}(T_{N}(W))/\widetilde{J}_{T_{N}(W)}$ and $W$ are irreducible, 
this surjective $V$-module map must be an equivalence. So we obtain 
$[S^{N}(T_{N}(W))/\widetilde{J}_{T_{N}(W)}]=[W]$, that is, 
 $g(f([W]))=[W]$. So we obtain $g\circ f=1_{[\mathfrak{W}]_{\rm irr}}$. 
This finishes the proof that $[\mathfrak{W}]_{\rm irr}$ is in bijection with $[\mathfrak{M}^{N}]_{\rm irr}$.
\epfv

\begin{cor}
For $N_{1}, N_{2}\in \N$ or equal to $\infty$, the set of  the equivalence classes of irreducible nondegenerate graded 
$A^{N_{1}}(V)$-modules is in bijection with the set of  the equivalence classes of irreducible 
nondegenerate
graded $A^{N_{2}}(V)$-modules. 
\end{cor}

We now assume that $V$ is a  M\"{o}bius vertex algebra, that is, 
a grading-restricted vertex algebra equipped with an operator $L_{V}(1)$ such that 
$L_{V}(1)$, $L_{V}(0)$ and $L_{V}(-1)$ satisfying the usually commutator relations for the 
standard basis of $\mathfrak{sl}_{2}$ and the usual commutator formula between
$L_{V}(1)$ and vertex operators for a vertex operator algebra. See, for example, 
Definition 7.1 in \cite{H-exist-twisted-mod} for the precise definition. 
In this case, a lower-bounded generalized $V$-module should also have an operator $L_{W}(1)$ satisfying the 
same relations as $L_{V}(1)$. We assume that $V$ is a grading-restricted 
M\"{o}bius vertex algebra in the remaining part of the paper because in this case,
a lowest weight of a lower-bounded generalized 
$V$-module is well defined. See Remark 7.3 in \cite{H-exist-twisted-mod}. 

\begin{prop}\label{l-b-f-l=g-r}
Let $V$ be a  M\"{o}bius vertex algebra.  Assume that $A^{N}(V)$ for all $N\in \N$ are finite dimensional
(for example, when $V$ is $C_{2}$-cofinite and of positive energy by Theorem \ref{finite-d}). Then 
every irreducible lower-bounded generalized $V$-module is an ordinary $V$-module and every lower-bounded generalized 
$V$-module of finite length is grading restricted. 
\end{prop}
\pf
Since for $N\in \N$, $A^{N}(V)$ is finite dimensional, there are only finitely many irreducible $A^{N}(V)$-modules. 
By Theorem \ref{irred-m-bij}, there are also finitely many irreducible lower-bounded generalized $V$-modules. 
For an irreducible 
lower-bounded generalized $V$-module $W$ with lowest weight $h_{W}$
 and $N\in \N$,  $T_{N}(W)$  is an irreducible nondegenerate graded $A^{N}(V)$-module
by Proposition \ref{T-N-preserv-irr}. Since $A^{N}(V)$ is finite dimensional, 
$T_{N}(W)$ is also finite dimensional. Thus $G_{N}(W)=W_{[h_{W}+N]}\subset T_{N}(W)$ is also finite 
dimensional. Since this is true for $N\in \N$, we see that $W$ is grading restricted. 
Since $W$ is irreducible, $L_{W}(0)$ must act semisimply on $W$. So $W$ is an irreducible ordinary $V$-module. 

Since as a graded vector space, a  lower-bounded generalized 
$V$ module $W$ of finite length  is a finite sum of irreducible lower-bounded generalized $V$-modules, 
which are all ordinary $V$-modules from what we have proved above. Then $W$ must be grading restricted.
\epfv

Since $V$ is a M\"{o}bius verex algebra, the associative algebras $A^{\infty}(V)$ and $A^{N}(V)$ for $N\in \N$ 
have an additional operator $L_{V}(1)$ induced from the operator $L_{V}(1)$ acting on 
$V$. For a lower-bounded generalized $V$-module $W$, there is also an 
operator $L_{Gr(W)}(1)$ on the $A^{\infty}(V)$-module $Gr(W)$ induced from 
$L_{W}(1)$ on $W$ such that $L_{Gr(W)}(1)$ maps $Gr_{n}(W)$ to $Gr_{n-1}(W)$. 
Restricting $L_{Gr(W)}(1)$ to $Gr^{N}(W)$, we obtain an operator $L_{Gr^{N}(W)}(1)$ on $Gr^{N}(W)$.

\begin{defn}
{\rm Let $V$ be a M\"{o}bius vertex algebra. A {\it graded $A^{N}(V)$-module}
is a graded $A^{N}(V)$-module $M$ when $V$ is viewed as a grading-restricted vertex algebra together with 
an operator $L_{M}(1)$ satisfying the following conditions:
\begin{enumerate}

\item $L_{M}(1)$ 
maps $G_{n}(M)$ to $G_{n-1}(M)$ for $n=0, \dots, N$, where $G_{-1}(M)=0$.

\item The operators $L_{M}(1)$ satisfies the commutator relations
\begin{align*}
{[L_{M}(0), L_{M}(1)]}&=-L_{M}(1),\nn
{[L_{M}(1), L_{M}(-1)]}&=2L_{M}(0),\nn
{[L_{M}(1), \vartheta_{M}([v]_{kl}+Q^{\infty}(V))]}&=\vartheta_{M}([(L_{V}(1)+2L_{V}(0)
+L_{V}(-1))v]_{(k-1)l}+Q^{\infty}(V)).
\end{align*}

\end{enumerate}
An graded $A^{N}(V)$-module $M$ is said to be {\it nondegenerate} if $M$ is nondegenerate 
when $V$ is viewed as a grading-restricted vertex algebra.
Let $M_{1}$ and $M_{2}$ be graded $A^{N}(V)$-modules.
An {\it graded $A^{N}(V)$-module map} from 
$M_{1}$ to $M_{2}$ is an $A^{N}(V)$-module map $f: M_{1}\to M_{2}$ 
such that $f(G_{n}(M_{1}))\subset G_{n}(M)$ for $n=0 \dots, N$, 
$f\circ L_{M_{1}}(1)=L_{M_{2}}(1)\circ f$,
$f\circ L_{M_{1}}(0)=L_{M_{2}}(0)\circ f$ and $f\circ L_{M_{1}}(-1)=L_{M_{2}}(-1)\circ f$.
A {\it graded $A^{N}(V)$-submodule} of a graded $A^{N}(V)$-module $M$ is an $A^{N}(V)$-submodule 
$M_{0}$ of $M$ such that  with the $A^{N}(V)$-module structure and the $\N$-grading induced from $M$
and the operators $L_{M}(1)\mbar_{M_{0}}$, $L_{M}(0)\mbar_{M_{0}}$ and $L_{M}(-1)\mbar_{M_{0}}$, 
$M_{0}$ is a graded $A^{N}(V)$-module.
A graded $A^{\infty}(V)$-module $M$ is said to be {\it generated by a subset $S$}  if 
$M$ is equal to the smallest graded $A^{N}(V)$-submodule containing $S$, or equivalently, 
$M$  is spanned by homogeneous elements 
obtained by applying elements of $A^{N}(V)$, $L_{M}(1)$ and
 $L_{M}(-1)$ to homogeneous summands of elements of $S$.  {\it Irreducible} and 
{\it completely reducible} graded $A^{N}(V)$-module 
are defined in the same way as in the case that $V$ is a grading-restricted vertex algebra.}
\end{defn}

From Proposition \ref{A-N-V-mod} and the property of $L_{W}(1)$, we immediately obtain the following:

\begin{prop}
Let $V$ be a M\"{o}bius vertex algebra.
For a lower-bounded generalized $V$-module $W$, $Gr^{N}(W)$
is nondegenerate graded $A^{N}(V)$-module.  Let $W_{1}$ and $W_{2}$ be lower-bounded generalized $V$-modules 
and $f: W_{1}\to W_{2}$ a $V$-module map. Then $f$ induces a graded $A^{N}(V)$-module map
$Gr^{N}(f): Gr^{N}(W_{1})\to Gr^{N}(W_{2})$.
\end{prop}

As is mentioned above, 
in the remaining part of this paper, we assume that $V$ is  a M\"{o}bius vertex algebra. We shall not 
repeat this assumption except in the statements of propositions, theorems, corollaries and so on. 
Lower-bounded generalized $V$-modules and graded $A^{N}(V)$-modules always mean 
those for $V$ as a M\"{o}bius vertex algebra, not as a grading-restricted vertex algebra. 
All the results that we have obtained above certainly still hold. 

We recall the notion of  lower-bounded generalized $V$-module 
of finite length. A lower-bounded generalized $V$-module $W$  is said to be {\it of fnite length} if 
there is a composition series $W=W_{0}\supset\cdots\supset W_{l+1}=0$ of 
lower-bounded generalized $V$-modules such that $W_{i}/W_{i+1}$ for $i=0, \dots, l$ 
are irreducible lower-bounded generalized $V$-modules. 

\begin{prop}\label{generators}
Let $V$ be a  M\"{o}bius vertex algebra. Assume that
the differences between the real parts of the lowest weights of the irreducible lower-bounded generalized
$V$-modules are all less than or equal to $N\in \N$. Then a lower-bounded generalized 
$V$-module $W$ of finite length
is generated by 
$$\coprod_{\Re(h_{W})\le \Re(n)\le \Re(h_{W})+N}W_{[n]}\subset \Omega_{N}(W),$$
where $h_{W}$ is a lowest weight of $W$.
\end{prop}
\pf
Let $W=W_{0}\supset W_{1}\supset \cdots \supset W_{l+1}=0$ be a finite composition series 
such that $W_{i}/W_{i+1}$ for $i=0, \dots, l$ are irreducible lower-bounded generalized $V$-modules. As a 
graded vector space, $W$ is isomorphic to $\coprod_{i=0}^{l}W_{i}/W_{i+1}$.
 In particular, the lowest weight of one of the irreducible lower-bounded generalized $V$-modules
$W_{i}/W_{i+1}$ for $i=0, \dots, l$ is a lowest weight $h_{W}$ of $W$.  

Let $w_{i}\in W_{i}$ be homogeneous for $i=0, \dots, l$ such that $w_{i}+W_{i+1}$ is 
a lowest weight vector of $W_{i}/W_{i+1}$. Then by assumption, the differences between the real parts of the lowest 
weights of $W_{i}/W_{i+1}$ for $i=0, \dots, l$ are less than or equal to $N$. Since one of these 
lowest weights is a lowest weight $h_{W}$ of $W$, we see that the the differences between the
real parts of the lowest weights of $W_{i}/W_{i+1}$ for $i=0, \dots, l$ and $\Re(h_{W})$
 are less than or equal to $N$. 
In particular $w_{i}\in \coprod_{\Re(h_{W})\le \Re(n)\le \Re(h_{W})+N}W_{[n]}$. 
Since for each $i$, $W_{i}/W_{i+1}$ is generated by 
$w_{i}+W_{i+1}$, $W_{i}$ is generated by $w_{i}$ and $W_{i+1}$. Thus 
$W$ is generated by $w_{i}$ for $i=0, \dots, l$. Since  $w_{i}\in \coprod_{\Re(h_{W})\le \Re(n)\le \Re(h_{W})+N}W_{[n]}$,
$W$ is generated by $\coprod_{\Re(h_{W})\le \Re(n)\le \Re(h_{W})+N}W_{[n]}$. It is clear that 
$\coprod_{\Re(h_{W})\le \Re(n)\le \Re(h_{W})+N}W_{[n]}$ is a subspace of $\Omega_{N}(W)$.
\epfv


\begin{cor}\label{cor-generators}
Let $V$ be a  M\"{o}bius vertex algebra.  Assume that $A_{N'}(V)$ for all $N'\in \N$ are finite dimensional
(for example, when $V$ is $C_{2}$-cofinite and of positive energy by Theorem \ref{finite-d}). Let $N\in \N$ such that 
the differences between the real parts of the lowest weights of the finitely many (inequivalent) 
irreducible ordinary $V$-modules are less than or equal to $N$.  Then a lower-bounded 
generalized $V$-module $W$ of finite length or a grading-restricted generalized $V$-module $W$
is generated by 
$$\coprod_{\Re(h_{W})\le \Re(n)\le \Re(h_{W})+N}W_{[n]}\subset \Omega_{N}(W).$$
\end{cor}
\pf
Since by Proposition \ref{l-b-f-l=g-r}, the finitely many (inequivalent)  irreducible lower-bounded generalized 
$V$-modules are all ordinary $V$-modules, the condition in Proposition 
\ref{generators} is satisfied. Also, by Corollary 3.16 in \cite{H-cofiniteness}, every grading-restricted 
generalized $V$-module is of finite length. Thus $W$ 
is generated by $\coprod_{\Re(h_{W})\le \Re(n)\le \Re(h_{W})+N}W_{[n]}$.
\epfv

\begin{thm}\label{comp-red-Gr-N}
Let $V$ be a M\"{o}bius vertex algebra. Assume that
the differences between the real parts of the lowest weights of the irreducible lower-bounded generalized
$V$-modules are all less than or equal to $N\in \N$. Then 
a lower-bounded generalized $V$-module $W$ of finite length is irreducible or completely reducible
 if and only if the  nondegenerate  graded $A^{N}(V)$-module $Gr^{N}(W)$  is irreducible or completely reducible, respectively.
\end{thm}
\pf
By Proposition  \ref{T-N-preserv-irr}, we already know that if $W$ is irreducible, $Gr^{N}(W)=T_{N}(W)$ is irreducibile. 
Conversely, 
assume that the nondegenerate  graded $A^{N}(V)$-module $Gr^{N}(W)$ is 
irreducible. Let $W_{0}$ be  a nonzero generalized $V$-submodule of $W$. 
Let $e_{W_{0}}: W_{0}\to W$ be the embedding map. Then we have a graded $A^{N}(V)$-moudle 
map $Gr(e_{W_{0}}): Gr^{N}(W_{0})\to Gr^{N}(W)$ given by 
$(Gr(e_{W_{0}}))(w_{0}+\Omega_{n-1}(W_{0}))=
w_{0}+\Omega_{n-1}(W)$ for $n=0, \dots, N$ and $w_{0}\in \Omega_{n}(W_{0})$. 
Since $e_{W}$ is injective, $Gr(e_{W_{0}})$ is also injective . So  $(Gr(e_{W_{0}}))(Gr^{N}(W_{0}))$ 
is a nondegenerate  graded $A^{N}(V)$-submodule of $Gr^{N}(W)$. Since $W_{0}$ is nonzero, $Gr^{N}(W_{0})$ is nonzero.
Since $Gr^{N}(W)$ is irreducible and $Gr(e_{W_{0}})$ is injective,  
$(Gr(e_{W_{0}}))(Gr^{N}(W_{0}))$ is equal to 
$Gr^{N}(W)$.  We now prove $W_{0}=W$. 
In fact, for $n=0, \dots, N$,  $(Gr(e_{W_{0}}))(Gr_{n}(W_{0}))
=\{w_{0}+\Omega_{n}(W)\mid w_{0}\in \Omega_{n}(W_{0})\}$. So 
$Gr_{n}(W)=\{w_{0}+\Omega_{n-1}(W)\mid w_{0}\in \Omega_{n}(W_{0})\}$.
For $n=0$, we obtain $\Omega_{0}(W)=Gr_{0}(W)=Gr(W_{0})=\Omega_{0}(W_{0})$. Assume that 
$\Omega_{n-1}(W)=\Omega_{n-1}(W_{0})$ for $n<N$. Given $w\in \Omega_{n}(W)$,
$w+\Omega_{n-1}(W)\in Gr_{n}(W)$. By
$Gr_{n}(W)=\{w_{0}+\Omega_{n-1}(W)\mid w_{0}\in \Omega_{n}(W_{0})\}$,
there exists $w_{0}\in \Omega_{n}(W_{0})$ such that $w+\Omega_{n-1}(W)=w_{0}+\Omega_{n-1}(W)$,
or equivalently, $w-w_{0}\in \Omega_{n-1}(W)=\Omega_{n-1}(W_{0})$. 
Thus $w\in \Omega_{n}(W_{0})$. This shows $\Omega_{n}(W)=\Omega_{n}(W_{0})$ for $n=0, \dots, N$.
In particular, $\Omega_{N}(W)=\Omega_{N}(W_{0})$. But by Proposition \ref{generators}, 
$W_{0}$ and $W$ are generated by $\Omega_{N}(W_{0})$ and $\Omega_{N}(W)$, respectively. 
Since $\Omega_{N}(W)=\Omega_{N}(W_{0})$, we must have $W=W_{0}$. 
So $W$ is irreducible. 

If $W$ is completely reducible, then by Proposition  \ref{T-N-preserv-irr}, 
$Gr^{N}(W)=T_{N}(W)$ is completely reducible. 
Conversely, assume that the nondegenerate  graded $A^{N}(V)$-module 
$Gr^{N}(W)$ is completely reducible. Then $Gr^{N}(W)=\coprod_{\mu\in \mathcal{M}}M^{\mu}$,
where $M^{\mu}$ for $\mu\in \mathcal{M}$ are irreducible nondegenerate  
graded $A^{N}(V)$-submodules of $Gr^{N}(W)$.
For $\mu\in \mathcal{M}$, since $M^{\mu}$ is  a nondegenerate  
graded $A^{N}(V)$-submodule of $Gr^{N}(W)$,
we have $M^{\mu}_{n}\subset 
Gr_{n}(W)=\Omega_{n}(W)/\Omega_{n-1}(W)$ for $n=0, \dots, N$. Let $W^{\mu}$ be the generalized $V$-submodule 
of $W$ generated by the set of 
elements of the form $w^{\mu}\in \Omega_{n}(W)$ such that $w^{\mu}+\Omega_{n-1}(W)\in M^{\mu}_{n}$
for  for $n=0, \dots, N$. Since $W^{\mu}$ is a generalized $V$-submodule of $W$, 
for $v\in V$, $k, l\in \N$ and $w^{\mu}\in \Omega_{l}(W)$ such that $w^{\mu}+\Omega_{l-1}(W)\in M^{\mu}_{l}$,
$$\res_{x}x^{l-k-1}Y_{W}(x^{L_{V}(0)}v, x)w^{\mu}+\Omega_{k-1}(W)\in Gr_{k}(W^{\mu}).$$
By the definition of $W^{\mu}$, we see that 
$\res_{x}x^{l-k-1}Y_{W}(x^{L_{V}(0)}v, x)w^{\mu}\in W^{\mu}$.
Since $w^{\mu}\in \Omega_{l}(W)$, $\res_{x}x^{l-k-1}Y_{W}(x^{L_{V}(0)}v, x)w^{\mu}=0$
for $k\in -\Z_{+}$. Therefore
$\res_{x}x^{l-k-1}Y_{W}(x^{L_{V}(0)}v, x)w^{\mu}\in W^{\mu}$ for $k\in \N$ are all the nonzero 
coefficients of $Y_{W}(v, x)w^{\mu}$.  So $W^{\mu}$ is closed under the action of 
the vertex operators on $W$. Since $M^{\mu}$ is invariant under the actions of $L_{Gr(W)}(0)$ and $L_{Gr(W)}(-1)$ 
and is a direct sum of generalized eigenspaces of $L_{Gr(W)}(0)$,
$W^{\mu}$ is invariant under the actions of $L_{W}(0)$ and $L_{W}(-1)$ 
and is a direct sum of generalized eigenspaces of $L_{W}(0)$. Thus $W^{\mu}$ is a generalized $V$-submodule of $W$. 

Let $w^{\mu}+\Omega_{n-1}(W^{\mu})\in Gr_{n}(W^{\mu})$,  where $0\le n\le N$ and 
$w^{\mu}\in \Omega_{n}(W^{\mu})
\subset \Omega_{n}(W)$. By the definition of $W^{\mu}$, we see that since $w^{\mu}$ is
 an element of $W^{\mu}$, 
$w^{\mu}+\Omega_{n-1}(W)\in G_{n}(M^{\mu})$. So we obtain a linear map from 
$Gr_{n}(W^{\mu})$ to $G_{n}(M^{\mu})$ given by $w^{\mu}+\Omega_{n-1}(W^{\mu})\mapsto
w^{\mu}+\Omega_{n-1}(W)$ for $w^{\mu}+\Omega_{n-1}(W^{\mu})\in Gr_{n}(W^{\mu})$. 
These maps for $n=0, \dots, N$ give a map from $Gr^{N}(W^{\mu})$ to $M^{\mu}$. 
It is clear that this map is a graded $A^{N}(V)$-module map. 
If for $0\le n\le N$,
the image $w^{\mu}+\Omega_{n-1}(W)$ of $w^{\mu}+\Omega_{n-1}(W^{\mu})\in Gr_{n}(W^{\mu})$
under this map is $0$ in $M^{\mu}$, then $w^{\mu}\in \Omega_{n-1}(W)$.
But $w^{\mu}\in \Omega_{n}(W^{\mu})\subset W^{\mu}$. So $w^{\mu}\in \Omega_{n-1}(W^{\mu})$
and $w^{\mu}+\Omega_{n-1}(W^{\mu})$ is $0$ in $Gr^{N}(W^{\mu})$. This means that this 
graded $A^{N}(V)$-module map is injective. In particular, the image of $Gr^{N}(W^{\mu})$ under this map 
is a nonzero nondegenerate  graded $A^{N}(V)$-submodule of $M^{\mu}$. But $M^{\mu}$ is irreducible.
So $Gr^{N}(W^{\mu})$ must be equivalent to $M^{\mu}$ and is therefore also irreducible. 
From what we have proved above, since $Gr^{N}(W^{\mu})$ is irreducible, 
$W^{\mu}$ is irreducible. Thus $W$ is complete reducible. 
\epfv

From Corollary \ref{cor-generators} and 
Theorem \ref{comp-red-Gr-N}, we obtain the following result:

\begin{cor}\label{rat-Gr-N}
Let $V$ be a  M\"{o}bius vertex algebra.  
Assume that $A_{N'}(V)$ for all $N'\in \N$ are finite dimensional
(for example, when $V$ is $C_{2}$-cofinite and of positive energy by Theorem \ref{finite-d}). Let $N\in \N$ such that 
the differences between the real parts of the lowest weights of the finitely many (inequivalent)
irreducible ordinary $V$-modules are less than or equal to $N$. 
Then a lower-bounded generalized $V$-module $W$ of finite length or a grading-restricted 
generalized $V$-module $W$ is a direct sum of irreducible ordinary $V$-modules
 if and only if the nondegenerate  graded $A^{N}(V)$-module $Gr^{N}(W)$  is completely reducible.
\end{cor}
\pf
Since $A^{N'}(V)$ is finite dimensional, there are only finitely many  (inequivalent)
 irreducible nondegenerate  graded $A^{N'}(V)$-modules. 
By Theorem \ref{irred-m-bij}, there are finitely many irreducible lower-bounded generalized  $V$-modules. 
By  Corollary \ref{cor-generators}, these finitely many irreducible lower-bounded generalized  $V$-modules
are all irreducible ordinary $V$-modules. 
There exists $N\in \N$ such that 
the differences between the real parts of the lowest weights of the finitely many
irreducible ordinary $V$-modules are less than or equal to $N$. For such $N$, the condition in 
Theorem \ref{comp-red-Gr-N} holds. 
So by Theorem \ref{comp-red-Gr-N}, a lower-bounded generalized $V$-module $W$ of 
finite length is a direct sum of irreducible ordinary $V$-modules
 if and only if $Gr^{N}(W)$  is completely reducible as a nondegenerate  graded
$A^{N}(V)$-module.

By Corollary 3.16 in \cite{H-cofiniteness}, every grading-restricted generalized $V$-module 
is of finite length. Thus the conclusion holds also for a grading-restricted 
generalized $V$-module $W$.
\epfv

\begin{rema}
{\rm Note that the assumption or condition on the lowest weights of irreducible 
$V$-modules in Proposition \ref{generators}, Corollary 
\ref{cor-generators}, Theorem \ref{comp-red-Gr-N} and Corollary \ref{rat-Gr-N} can be 
weakened to the assumption or condition that the differences between the real parts of the 
lowest weights of the irreducible lower-bounded generalized
$V$-modules appearing as a quotient in a composition series of 
$W$ are all less than or equal to $N\in \N$. This is because the proofs used only this weaker assumption
or condition. For the study of some particular 
lower-bounded generalized $V$-modules of finite length or some grading-restricted 
generalized $V$-modules, this weaker assumption or condition is certainly easier to verify 
than the more general ones in the statements of these results.}
\end{rema}

\vspace{1em}

\noindent {\small \sc Department of Mathematics, Rutgers University,
110 Frelinghuysen Rd., Piscataway, NJ 08854-8019}

\noindent {\em E-mail address}: {\tt yzhuang@math.rutgers.edu}

\end{document}